\documentclass{article}

\usepackage{latexsym, amssymb}

\newtheorem{definition}{Definition}[section]
\newtheorem{theorem}[definition]{Theorem}

\newtheorem{notation}[definition]{Notation}
\newtheorem{lemma}[definition]{Lemma}
\newtheorem{conjecture}[definition]{Conjecture}
\newtheorem{remark}[definition]{Remark}
\newtheorem{corollary}[definition]{Corollary}
\newtheorem{problem}[definition]{Problem}
\def\R{\mathbb R}
\def\C{\mathbb C}

\newcommand{\beast}{\begin{eqnarray*}}
\newcommand{\eeast}{\end{eqnarray*}}

\begin{document}

\title{Taut distance-regular graphs and the subconstituent algebra}
\author{Mark S. MacLean and Paul Terwilliger}

\maketitle


\begin{abstract}
We consider a bipartite distance-regular graph $\Gamma $
with diameter $D \geq 4$ and valency $k \ge 3$. 
Let
$X$ denote the vertex set of $\Gamma $
and fix $x \in X$.  Let $\Gamma_{2}^{2}$ denote the
graph with vertex set $\breve{X}= \{y \in X\;|\;\partial(x,y)=2\},$ and 
edge set $\breve{R}= \{yz\;|\;y,z \in \breve{X},\, \partial(y,z)=2 \},$ 
where $\partial$ is the path-length distance function for $\Gamma$.
The graph $\Gamma_{2}^{2}$ has 
exactly $k_{2}$ vertices, where $k_{2}$ is the second valency of
$\Gamma.$  
Let
$\eta_1, \eta_2 , \ldots , \eta_{k_{2}}$ denote the eigenvalues of the 
adjacency matrix of $\Gamma_{2}^{2}$; we call these the {\it local 
eigenvalues of $\Gamma$.}  Let $A$ denote the adjacency matrix of 
$\Gamma$.  We obtain upper and lower bounds for the 
local eigenvalues in terms of the intersection numbers of $\Gamma$ and the eigenvalues 
of $A$.  
Let
$T=T(x)$
denote the subalgebra
of
$\hbox{Mat}_X(\C)$ generated by $A, E^*_0, E^*_1, \ldots, E^*_D$, 
where for $0 \le i \le D$, $E^*_i$ represents the 
projection onto the $i^{\hbox{th}}$ subconstituent of $\Gamma $ with 
respect to $x$.
We refer to $T$ as the subconstituent algebra (or Terwilliger algebra) of 
$\Gamma$ with respect to $x$.
An irreducible $T$-module $W$ is said to be {\it thin}
whenever $\hbox{dim}E^*_iW\leq 1$ for
$0 \leq i \leq D$.
By the {\it endpoint} of $W$ we mean $\hbox{min}\lbrace i 
|E^*_iW\not=0 \rbrace $.  We give a detailed description of the 
thin irreducible $T$-modules that have 
endpoint 2 and dimension $D-3$.  
In [{\em Discrete Math.}, 225(2000), 193--216] MacLean defined what it means for $\Gamma$ to be 
{\it taut}.  We obtain three characterizations of the 
taut condition, each of which involves the local eigenvalues or the 
above $T$-modules.
\end{abstract}



\section{Introduction}

\bigskip \noindent
Let $\Gamma$ denote a distance-regular graph with diameter $D \geq 
4$, valency $k \ge 3$, and intersection numbers $a_i, b_i, c_i$ 
(see Section 2 for formal definitions).
In this paper we obtain some results on the subconstituent
algebra \cite{terwSub1} of $\Gamma$,
and some related results concerning the taut condition \cite{maclean1}.
We will state our results shortly, but first we motivate
things with a brief discussion of the subconstituent algebra.
Let
$X$ denote the vertex set of $\Gamma$ and fix $x 
\in X$. We view $x$ as a ``base vertex.'' Let $T=T(x)$ denote the subalgebra of
$\hbox{Mat}_X(\C)$ generated by $A, E^*_0, E^*_1, \ldots, E^*_D$, 
where
$A$ denotes the adjacency matrix of $\Gamma $ and $E^*_i$ represents the 
projection onto the $i^{\hbox{th}}$ subconstituent of $\Gamma $ with 
respect to $x$. The algebra $T$ is called the {\it subconstituent 
algebra} (or {\it Terwilliger algebra}) of $\Gamma $
with respect to $x$
\cite{terwSub1}. Observe $T$ has finite dimension.
Moreover $T$ is semi-simple; the reason is
each of $A, E^*_0, E^*_1, \ldots, E^*_D$ is symmetric with
real entries, so
$T$ is closed under the conjugate-transpose map \cite[p. 157]{CR}. 
Since $T$ is semi-simple, each $T$-module is a direct sum of 
irreducible $T$-modules. Describing the irreducible $T$-modules is an
active area of research
\cite{caugh2}--\cite{curtin6}, \cite{dickie1}--\cite{gotight}, 
\cite{hobart},
\cite{tanabe},
\cite{terwSub1}--\cite{terwSubEndpt1},
\cite{tomiyama2}.

\bigskip \noindent
In this paper we are concerned with the irreducible $T$-modules that 
possess a certain property.  In order to define this property we make 
a few observations.
Let $W$ 
denote an irreducible $T$-module.
Then $W$ is the direct sum of the
nonzero spaces among $E^*_0W, E^*_1W,\ldots, E^*_DW$. There is a 
second decomposition of interest.
To obtain it we make a definition.
Let $k=\theta_0 > \theta_1 > \cdots > \theta_D$
denote the distinct eigenvalues of $A$, and for $0 \leq i \leq D $ 
let $E_i$
denote the primitive idempotent of $A$ associated with $\theta_i$.
Then $W$ is the direct sum of the
nonzero spaces among $E_0W, E_1W,\ldots, E_DW$.
If the dimension of $E^*_iW$ is at most 1 for
$0 \leq i \leq D$ then
the dimension of $E_iW$ is at most 1 for
$0 \leq i \leq D$ \cite[Lemma 3.9]{terwSub1};
in this case we say $W$ is {\it thin}.
Let $W$ denote an irreducible $T$-module.
By the {\it endpoint} of $W$
we mean $\hbox{min}\lbrace i | 0 \leq i \leq D, \;E^*_iW\not=0 
\rbrace $.
There exists a unique irreducible $T$-module with endpoint 0 
\cite[Proposition 8.4]{egge1}. We call this module $V_0$.
The module $V_0$ is thin; in fact
$E^*_iV_0$ and
$E_iV_0$ have dimension 1 for $0 \leq i \leq D$ \cite[Lemma 
3.6]{terwSub1}. For a detailed description of $V_0$ see
\cite{curtin1}, \cite{egge1}. 

\bigskip \noindent
For the rest of this section assume $\Gamma$ is bipartite.  
There exists, up to isomorphism, a unique irreducible $T$-module with endpoint 1 
\cite[Corollary 7.7]{curtin1}.  We call this module $V_{1}$.  
The module $V_1$ is thin; in fact each of
$E^*_iV_1$,
$E_iV_1$ has dimension 1 for $1 \leq i \leq D-1$ and 
$E^{*}_{D}V_{1}=0$, $E_{0}V_{1}=0$, $E_{D}V_{1}=0$. 
For a detailed description of $V_1$ see
\cite{curtin1}. In this paper we are concerned with the 
thin irreducible $T$-modules with endpoint 2.

\bigskip \noindent
We now state our results.
In order to state our first result
we define some parameters.
Let $\Gamma_{2}^{2} = \Gamma_{2}^{2}(x) $ denote the
graph with vertex set $\breve{X}$ and edge set $\breve{R}$, where
\begin{eqnarray}
\breve{X} &=& \{y \in X\;|\;\partial(x,y)=2\},\nonumber\\
\breve{R} &=& \{yz\;|\;y,z \in \breve{X},\, \partial(y,z)=2 \},\nonumber
\end{eqnarray}
and where $\partial$ is the path-length distance function for $\Gamma$.
The graph $\Gamma_{2}^{2}$ has 
exactly $k_{2}$ vertices, where $k_{2}$ is the second valency of
$\Gamma.$ Also, $\Gamma_{2}^{2}$ is regular with valency $p^{2}_{22}$.
We let
$\eta_1, \eta_2 , \ldots , \eta_{k_{2}}$
denote the eigenvalues of the adjacency matrix of $\Gamma_{2}^{2}$. 
By \cite[Theorem 11.7]{curtin2}, these eigenvalues may be ordered such that $\eta_{1} = p_{22}^{2}$ 
and $\eta_{i}= b_{3}-1$ $(2 \leq i \leq k)$.   

\bigskip \noindent
Abbreviate $d = \lfloor D/2 \rfloor$.  Our first main result is that 
${\tilde \theta}_1 \leq 
\eta_i \leq {\tilde \theta}_d $ for $k+1 \leq i \leq k_{2}$, where
${\tilde \theta}_1=-1-b_2 b_{3}(\theta_1^{2}-b_{2})^{-1}$ and 
${\tilde \theta}_d=-1-b_2 b_{3}(\theta_d^{2}-b_{2})^{-1}$. 
We remark 
$\theta_{1}^{2} > b_{2} > \theta_{d}^{2}$ by \cite[Lemma 2.6]{jkt},
so ${\tilde \theta}_1 < -1$ and ${\tilde \theta}_d\geq 0$.

\bigskip \noindent
In order to state our next set of results we make a definition. 
Let $W$ denote a thin irreducible $T$-module with endpoint 2.
Observe $E^*_2W$ is a $1$-dimensional eigenspace for $E^*_2A_{2}E^*_2$;
let $\eta $ denote the corresponding eigenvalue.
It turns out $\eta $
is one of $\eta_{k+1}, \eta_{k+2},\ldots, \eta_{k_{2}}$ so
${\tilde \theta}_1 \leq \eta \leq {\tilde \theta}_d$.
We call $\eta $ the {\it local eigenvalue } of $W$.
To describe the structure of $W$ we distinguish four cases: (i) $D$ 
is odd, and $\eta 
={\tilde \theta}_1 $ or $\eta ={\tilde \theta}_d $; (ii) $D$ is even 
and $\eta = {\tilde \theta}_{1}$; (iii) $D$ is even 
and $\eta = {\tilde \theta}_{d}$; (iv) $ {\tilde \theta}_1 < \eta < {\tilde 
\theta}_d$.  
We investigate cases (i), (ii) in the present paper.
We will investigate the remaining cases in a future paper.

\bigskip \noindent
Our results concerning the $T$-modules are as follows. 
Choose $n \in \{1,d\}$ if $D$ is odd, and let $n=1$ if $D$ is even.  
Define $\eta = {\tilde \theta}_{n}$.  
Let $W$ denote a thin irreducible
$T$-module with endpoint 2 and local eigenvalue $\eta$.  
Let $v$ denote a nonzero vector in $E^*_2W$.
We show $W$ has a basis
$E_iv$ $(1 \leq i \leq D-1,\;i\not=n,\; i \not= D-n)$.
We show this basis is orthogonal (with respect to the Hermitian dot
product) and we compute the square norm of each basis vector.
We show $W$ has a basis
$E^*_{i+2}A_iv$ $(0 \leq i \leq D-4)$,
where $A_i$ denotes the $i^{\hbox{th}}$
distance matrix for $\Gamma $.
We find the matrix representing $A$ with respect to this basis.
We show this basis is orthogonal and we compute the square norm of 
each basis vector. We find the transition matrix relating our two 
bases for $W$.
We show the following scalars are equal: (i) The multiplicity with 
which $W$ appears
in the standard module $\mathbb{C}^{X}$; (ii) The number of times $\eta$ 
appears among
$\eta_{k+1}, \eta_{k+2}, \ldots, \eta_{k_{2}}$.

\bigskip \noindent
In order to state our remaining results we recall the taut condition.
In \cite[Theorem 18]{curtin3} Curtin showed that $b_{2}(k-2) \ge (c_{2}-1)  \theta_{1}^{2}
$ with equality if and only if $\Gamma$ is 
2-homogeneous in the sense of Nomura 
\cite{Nom1}.  In \cite[Theorem 12]{curtin3} Curtin showed $\Delta \ge 0$, where 
$
\Delta = (k-2)(c_{3}-1)-(c_{2}-1)p_{22}^{2}.
$
In \cite[Lemma 3.8]{maclean1} MacLean proved
that 
\begin{equation}\label{bfb2}
b_{3}\left(b_{2}(k-2)-(c_{2}-1)\theta_1^{2}\right)\left( 
b_{2} (k-2)-(c_{2}-1)\theta_d^{2}\right) \ge 
b_{1}\Delta(\theta_{1}^{2}-b_{2})(b_{2}-\theta_{d}^{2}).
\end{equation}
We mentioned earlier that $\theta_{1}^{2} > b_{2} > \theta_{d}^{2}$, 
so the last two factors on the right in (\ref{bfb2}) are positive.
Observe each factor in (\ref{bfb2}) is nonnegative.    
From these comments we find that $\Gamma$ is
2-homogeneous if and only if $\Delta =0$ and equality holds in 
(\ref{bfb2}).  
MacLean defined $\Gamma$ to be {\it taut} \cite{maclean1} whenever $\Delta \not= 0$ 
and equality holds in (\ref{bfb2}). 

\bigskip \noindent
Assume for the moment that $\Gamma$ is taut.
It turns out that the structure of $\Gamma$ depends to a large extent on
the parity of $D$.
We investigated this structure for $D$ even in \cite{maclean3},
and for $D$ odd in \cite{maclean2}.  
In this paper we obtain three characterizations of the taut condition,
two of which require the assumption that $D$ is odd.

\bigskip \noindent
In order to state our first characterization of the taut condition 
we make a definition.  We say $\Gamma$ is {\it spectrally taut 
with respect to $x$} whenever $\eta_i$ is one of
${\tilde \theta}_1,
{\tilde \theta}_d $ for $k+1 \leq i \leq k_{2}$.  We show the 
following are equivalent:  (i) $\Gamma$ is taut; (ii) $\Delta \not= 0$ 
and $\Gamma$ is 
spectrally taut with respect to each vertex; (iii) $\Delta \not= 0$ 
and $\Gamma$ is 
spectrally taut with respect to at least one vertex.

\bigskip \noindent
For the rest of this section assume $D$ is odd.
We obtain two additional characterizations of the taut 
condition.  In order to state the first one we make a definition.  
We say $\Gamma$ is {\it taut with respect to $x$}
whenever every irreducible $T$-module with endpoint 2 is
thin with local eigenvalue
${\tilde \theta}_1$ or
${\tilde \theta}_d$.  
We show
the following are equivalent: (i) $\Gamma$ is taut;
(ii) $\Delta \not= 0$ and $\Gamma$ is taut with respect to each vertex;
(iii) $\Delta \not= 0$ and $\Gamma$ is taut with respect to at least one vertex.

\bigskip \noindent
Before we state our last characterization we recall two concepts.  
Recall that 
$\Gamma $ is {\it 2-thin with respect to $x$}
whenever every irreducible $T$-module with endpoint 2 is
thin.  Recall that $\Gamma$ is an {\it antipodal 2-cover} 
whenever 
for all $y \in X$, there exists a unique vertex $z \in X$ such that 
$\partial (y,z)=D$.  
We show the following are equivalent:
(i) $\Gamma$ is taut or 2-homogeneous;
(ii) $\Gamma$ is an antipodal 2-cover and 2-thin with respect to each vertex;
(iii) $\Gamma$ is an antipodal 2-cover and 2-thin with respect to at least one 
vertex. 

\bigskip \noindent
For more information on the taut condition and related 
topics we refer the reader to \cite{gotight}, \cite{jk}--\cite{aap4},
\cite{terwSubEndpt1},
\cite{tomiyama}.


\section{Preliminaries concerning distance-regular graphs}
In this section we review some definitions and basic concepts 
concerning distance-regular
graphs.
For more background information we refer the reader to \cite{bannai}, 
\cite{bcn}, \cite{godsil} or \cite{terwSub1}.

\bigskip
\noindent
Let $X$ denote a nonempty finite set.
Let
$\hbox{Mat}_X(\C)$ denote the $\C$-algebra
consisting of all matrices whose rows and columns are indexed by $X$
and whose entries are in $\;\C $. Let
$V=\C^X$ denote the vector space over $\C$
consisting of column vectors whose coordinates are indexed by $X$ and 
whose entries are
in $\C$.
We observe
$\hbox{Mat}_X(\C)$ acts on $V$ by left multiplication.
We endow $V$ with the Hermitian inner product $\langle \, , \, 
\rangle$ defined by
\begin{equation}\label{DOTDEF}
\langle u,v \rangle = u^t\overline{v} \qquad (u,v \in V),
\end{equation}
where $t$ denotes transpose and $-$ denotes complex conjugation.
As usual, we abbreviate $\|u\|^2 = \langle u,u \rangle$ for all $u 
\in V.$ For all $y \in X,$ let $\hat{y}$ denote the element
of $V$ with a 1 in the $y$ coordinate and 0 in all other coordinates.
We observe $\{\hat{y}\;|\;y \in X\}$ is an orthonormal basis for $V.$
The following formula will be useful. For all $B \in 
\hbox{Mat}_X(\C)$ and for all $u,v \in V$, \begin{equation}\label{ADJ}
\langle Bu,v \rangle = \langle u, {\overline B}^tv \rangle. 
\end{equation}

\bigskip
\noindent
Let $\Gamma = (X,R)$ denote a finite, undirected, connected graph,
without loops or multiple edges, with vertex set $X$ and edge set
$R$.  
Let $\partial $ denote the
path-length distance function for $\Gamma $, and set
$D = {\rm max}\,\{\partial(x,y) \;|\; x,y \in X\}. $ We refer to 
$D$ as the {\it diameter} of $\Gamma $.  We write $\lfloor D/2 
\rfloor$ to denote the greatest integer at most $D/2$.  
Let $x,y$ denote vertices of $\Gamma$. We say $x,y$ are {\it 
adjacent} whenever $xy$ is an edge.
Let $k$ denote a nonnegative integer. We say $\Gamma $ is
{\it regular} with {\it valency k} whenever each vertex of $\Gamma$ 
is adjacent to exactly $k$ distinct vertices of $\Gamma $.
We say $\Gamma$ is {\it distance-regular}
whenever for all integers $h,i,j\;(0 \le h,i,j \le D)$ and for all
vertices $x,y \in X$ with $\partial(x,y)=h,$ the number
\begin{equation}\label{PHIJ}
p_{ij}^h = |\{z \in X \; |\; \partial(x,z)=i, \partial(z,y)=j \}|
\end{equation}
is independent of $x$ and $y.$ The integers $p_{ij}^h$ are called
the {\it intersection numbers} of $\Gamma.$ We abbreviate
$c_i= p_{1i-1}^i \;(1 \le i \le D),\;a_i = p_{1i}^i \; (0 \le i \le 
D)$,
and $b_i= p_{1i+1}^i \; (0 \le i \le D-1)$. For notational
convenience we define $c_0=0$ and $b_D=0$.  We note $a_0=0$ and 
$c_1=1$.

\bigskip
\noindent
For the rest of this paper we assume $\Gamma$ is distance-regular 
with diameter $D\geq 3$.

\bigskip
\noindent
By (\ref{PHIJ}) and the triangle inequality, $
p_{ij}^h = 0$ if one of $h,i,j$ is bigger than the sum of the other 
two $(0 \le 
h,i,j \le D)$.
Observe $\Gamma$ is regular with valency $k=b_0,$ and that
$
c_i+a_i+b_i = k \quad (0 \le i \le D).$
Moreover $b_i > 0\;(0 \le i \le D-1)$ and $c_i>0 \;(1 \le i \le D).$

\bigskip
\noindent
For $0 \le i \le D$ we abbreviate $k_i=p_{ii}^0,$ and observe
\begin{equation}\label{DEFKI}
k_i = |\{z \in X \;|\; \partial(x,z)=i\}|,
\end{equation}
where $x$ is any vertex in $X$. Apparently $k_0=1$ and $k_1 =k.$
By \cite[p.195]{bannai} we have
\begin{equation}\label{KI}
k_i = \frac{b_0b_1\cdots b_{i-1}}{c_1c_2\cdots c_i} \qquad \qquad (0 
\le i \le D).
\end{equation}
We refer to $k_i$ as the {\it $i^{\hbox {th}}$ valency of } 
$\Gamma.$\\

\bigskip
\noindent
We recall the Bose-Mesner algebra of $\Gamma.$ For $0 \le i \le D$ 
let $A_i$ denote the matrix in $\hbox{Mat}_X(\C)$ with
$xy$ entry
$$
{(A_i)_{xy} = \cases{1, & if $\partial(x,y)=i$\cr
0, & if $\partial(x,y) \ne i$\cr}} \qquad (x,y \in X).
$$

\bigskip
\noindent
We call $A_i$ the $i^{\hbox{th}}$ {\it distance matrix} of $\Gamma.$
For convenience we define $A_i=0$ for $i < 0$ and $i > D.$
We abbreviate $A=A_1$ and call this the {\it adjacency
matrix} of $\Gamma.$ We observe (ai) $A_{0}=I$;  (aii) 
$\sum_{i=0}^{D}A_{i}=J$;  (aiii) $\overline{A}_i = A_i$ $(0 \le i \le 
D)$;  (aiv) $A_i^t = A_i$ $(0 \le i \le D)$;  (av) 
$A_iA_j = \sum_{h=0}^D p_{ij}^h A_h$ $( 0 \le i,j \le 
D)$, where $I$ denotes the identity matrix and $J$ denotes the all 1's 
matrix.
Let $M$ denote the subalgebra of $\hbox{Mat}_X(\C)$ generated by $A.$ 
Using (ai) and (av)
one can readily show $A_0,A_1,\ldots,A_D$
form a basis for $M.$ 
We refer to $M$ as the {\it Bose-Mesner algebra} of $\Gamma$.
By \cite[p. 45]{bcn} $M$ has a second basis 
$E_0,E_1,\ldots,E_D$ such that
(ei) $E_0 = |X|^{-1}J$;  (eii) 
$\sum_{i=0}^D E_i = I$;  (eiii)
$\overline{E}_i = E_i$ $(0 \le i \le D)$; (eiv)
$E_i^t = E_i$ $(0 \le i \le D)$;  (ev)
$E_iE_j = \delta_{ij}E_i$ $(0 \le i,j \le D)$.  
We refer to $E_0, E_1, \ldots, E_D $ as the {\it primitive 
idempotents}
of $\Gamma$. We call $E_0$ the {\it trivial idempotent} of $\Gamma.$

\bigskip
\noindent
We recall the eigenvalues
of $\Gamma $.
Since $E_0,E_1,\ldots,E_D$ form a basis for $M,$ there exist complex 
scalars $\theta_0,\theta_1,
\ldots,\theta_D$ such that
$A = \sum_{i=0}^D \theta_iE_i$.
Combining this with (ev)
we find $AE_i = E_iA = \theta_iE_i$ for $0 \leq i \leq D$.
Using (aiii) and
(eiii) we find $\theta_0,\theta_1,\ldots,\theta_D$ are in $\R.$ 
Observe
$\theta_0,\theta_1,\ldots,\theta_D$ are distinct since $A$ generates 
$M.$ By 
\cite[Proposition 3.1]{biggs} we have $\theta_0=k$
and $-k \le \theta_i \le k$ for $0 \le i \le D.$ Throughout this paper
we assume $E_0, E_1,\ldots, E_D$ are indexed so that $\theta_0 > 
\theta_1 >\cdots >\theta_D.$ We refer to $\theta_i$ as the {\it 
eigenvalue} of $\Gamma$ associated
with $E_i.$ 
For $0 \leq i \leq D$ let $m_i$ denote the rank of $E_i$. We refer to 
$m_i$ as the {\it multiplicity} of $E_i$ (or $\theta_i$).
From (ei) we find $m_0=1$. Using (eii)--(ev)
we find \begin{equation}\label{EIVDECOM}
V = E_0V+E_1V+ \cdots +E_DV \qquad \qquad {\rm (orthogonal\ direct\ 
sum}).
\end{equation}
For $0 \le i \le D$ the space $E_iV$ is the eigenspace of $A$ 
associated with $\theta_i$. We observe the dimension of $E_iV$ is 
$m_i$.

\bigskip
\noindent We now recall the dual eigenvalues of $\Gamma $.
Let $\theta$ denote an eigenvalue of $\Gamma$ and let $E$ denote
the associated primitive idempotent. Since $A_0,A_1,\ldots,A_D$
is a basis for $M,$ there exist complex scalars
$\theta_0^*,\theta_1^*, \ldots,\theta_D^*$ such that
\begin{equation}\label{PRIMIDEM}
E = |X|^{-1} \sum_{i=0}^D \theta_i^*A_i.
\end{equation}
Evaluating (\ref{PRIMIDEM}) using (aiii) and (eiii) we see 
$\theta_0^*,\theta_1^*, \ldots, \theta_D^*$ are in $\R.$
We refer to $\theta_i^*$ as the $i^{\hbox{th}}$ {\it dual eigenvalue} 
of $\Gamma$ with respect to $E$ (or $\,\theta$). We call 
$\theta_0^*,\theta_1^*, \ldots,\theta_D^*$ the {\it dual eigenvalue 
sequence} associated with $E$ (or $\,\theta$). By \cite[p. 128]{bcn} 
we have \begin{equation}\label{DRECUR}
c_i\theta_{i-1}^* + a_i\theta_i^* + b_i\theta_{i+1}^* = 
\theta\theta_i^* \qquad \qquad (0 \leq i \leq D),
\end{equation}
where $\theta_{-1}^*,\; \theta_{D+1}^*$ are indeterminates. We remark 
by \cite[p. 62]{bannai} that 
$\theta_0^* = m_i$ where
$\theta=\theta_i.$

\bigskip \noindent
The following lemma will be useful.

\begin{lemma}\cite[Lemma 2.6]{jkt}\label{THETAJ}
Let $\Gamma $ denote a distance-regular graph with diameter $D \ge 3$ 
and eigenvalues $k=\theta_0 > \theta_1 > \cdots >\theta_D.$ Then
\begin{description}
\item[(i)] $-1 < \theta_1 < k.$
\item[(ii)] $a_1-k \le \theta_D < -1.$
\end{description}
\end{lemma}

\noindent
Later in this paper we will discuss polynomials in one or two 
variables. We 
will use the following notation.
We let $\lambda $ denote an indeterminate. We let
$\R \lbrack \lambda \rbrack $ denote the $\R$-algebra consisting of 
all polynomials in $\lambda $ that
have coefficients in $\R$.  We let $\mu$ denote an indeterminate that 
commutes with $\lambda$.  We let $\R \lbrack \lambda, \mu \rbrack $
denote the $\R$-algebra consisting of 
all polynomials in $\lambda $ and $\mu$ that
have coefficients in $\R$.


\section{Bipartite distance-regular graphs}

\bigskip

        We now consider the case in which $\Gamma$ is bipartite.  
        We say $\Gamma$ is {\it bipartite}
        whenever the vertex set $X$ can be partitioned into two subsets, neither
        of which contains an edge. In the next few lemmas, we recall some routine
        facts concerning the case in which $\Gamma$ is bipartite.  
        To avoid trivialities, we will 
        generally assume $D \geq 4$.   
        
\begin{lemma} \label{bipequiv}  \cite[Propositions 3.2.3, 4.2.2]{bcn} \quad
        Let $\Gamma$ denote a distance-regular graph with diameter $D \geq 4$,
        valency $k$, 
        and eigenvalues
        $\theta_0 > \theta_1 > \cdots > \theta_D$.  The following are
        equivalent:
\begin{description}
\item[(i)] 
       $\Gamma$ is bipartite.
\item[(ii)]
       $p^{h}_{ij} = 0 \mbox{ if } h+i+j \mbox{ is odd} \qquad (0 \le 
       h,i,j \le D)$.
\item[(iii)]
        $a_i = 0 \qquad (0 \le i \le D)$.
\item[(iv)]
        $c_i + b_i =k  \qquad (0 \le i \le D)$.
\item[(v)]
        $\theta_{D-i} = -\theta_i \qquad (0 \le i \le D).$
\end{description} 
\end{lemma}

\begin{lemma}  \label{thetad}
Let $\Gamma$ denote a bipartite distance-regular graph with diameter $D \geq 4$
        and eigenvalues
        $k=\theta_0 > \theta_1 > \cdots > \theta_D$.
\begin{description}
\item[(i)] Assume $D$ is even and let $d = D/2$.  Then $\theta_{d}=0$.
\item[(ii)] Assume $D$ is odd and let $d =(D-1)/2$.  Then 
$\theta_{d} > 0$ and $\theta_{d+1}=-\theta_{d}$.
\end{description}
\end{lemma}        
{\it Proof.}  Immediate from Lemma \ref{bipequiv}(v).  \hfill $\Box$\\

\begin{lemma} \label{oppsigs}  \cite[Lemma 9]{curtin3}
        Let $\Gamma$ denote a bipartite distance-regular graph with
        diameter $D \geq 4$.  Let $\theta$ denote an eigenvalue of $\Gamma$ 
        and let $\theta^{*}_0, \theta^{*}_1, \ldots, 
        \theta^{*}_D$ denote the corresponding dual eigenvalue sequence.  Then the dual eigenvalue sequence associated
        with $-\theta$ is $\theta^{*}_0, -\theta^{*}_1, \theta^{*}_2, \ldots, 
        (-1)^D \theta^{*}_D$.  
\end{lemma}

\begin{lemma}  \label{J-ED}
   Let $\Gamma =(X,R)$ denote a bipartite distance-regular graph with diameter 
   $D \geq 4$,
        and eigenvalues
        $\theta_0 > \theta_1 > \cdots > \theta_D$.   Then $E_{D} = 
        |X|^{-1}J'$, where 
\begin{equation}  \label{J-def}
     J' = \sum_{i=0}^{D} (-1)^{i}A_{i}.
\end{equation}
\end{lemma}
{\it Proof.}  By (aii), (ei), we have $E_{0} = 
|X|^{-1}\sum_{i=0}^{D}A_{i}$.  Combining this with (\ref{PRIMIDEM}) 
and Lemma \ref{oppsigs} we find $E_{D} = 
|X|^{-1}\sum_{i=0}^{D}(-1)^{i}A_{i}$.  The 
result follows.
\hfill $\Box$\\

\begin{lemma}  \label{bipeigs}
Let $\Gamma$ denote a bipartite distance-regular graph with diameter 
$D \geq 4$ and eigenvalues $\theta_{0} > \theta_{1} > \cdots > 
\theta_{D}$.  Then $\theta_{1}^{2} > 
b_{2} > \theta_{d}^{2},$ where $d = \lfloor D/2 \rfloor$.
\end{lemma}
{\it Proof.}  Apply Lemma \ref{THETAJ} to the halved graph of 
$\Gamma$, and use \cite[Proposition 4.2.3]{bcn}.
\hfill $\Box$\\

\begin{lemma} \label{p222lem} \cite[Lemma 4.1.7]{bcn}
Let $\Gamma$ denote a bipartite distance-regular graph with diameter 
$D \geq 4$.  Then
\begin{equation} \label{p222}
   p^{2}_{22} = (b_{2}(c_{3}-1)+c_{2}(k-2)){c_{2}}^{-1}.
\end{equation}
\end{lemma}


\section{Two families of polynomials}

\bigskip
\noindent
Let $\Gamma=(X,R)$ denote a bipartite distance-regular graph with diameter $D 
\ge 4.$ In this section we recall two types of polynomials associated 
with
$\Gamma $. To motivate things, we recall by 
(av) and the triangle inequality that \begin{equation}\label{AAI1}
AA_i = b_{i-1}A_{i-1}  + c_{i+1}A_{i+1} \qquad \qquad
(0 \le i \le D),
\end{equation}
where $b_{-1}=0$ and $c_{D+1}=0$.
Let $f_0,f_1,\ldots,$$f_D$ denote the polynomials in $\R[\lambda]$ 
satisfying $f_0=1$ and
\begin{equation}\label{FIPOLY2}
\lambda f_i = b_{i-1}f_{i-1}  + c_{i+1}f_{i+1} \qquad (0 \le 
i \le D-1),
\end{equation}
where $f_{-1}=0$.
Let $i$ denote an integer $(0 \leq i \leq D)$. The polynomial $f_i$ 
has degree $i$, and
the coefficient of
$\lambda^i$ is $(c_1c_2\cdots c_i)^{-1}.$
Comparing (\ref{AAI1}) and (\ref{FIPOLY2}) we find
$f_i(A)=A_i$. By \cite[p. 63]{bannai}
the polynomials $f_0,f_1,\ldots, f_D$ satisfy the orthogonality 
relation
\begin{equation}\label{FIFJ1}
\sum_{h=0}^D f_i(\theta_h)f_j(\theta_h)m_h = \delta_{ij} |X|k_i
\qquad \qquad (0 \le i,j \le D).
\end{equation}

\bigskip
\noindent
Let $\theta$ denote an eigenvalue of $\Gamma$ and let
$\theta_0^*,\theta_1^*, \ldots,\theta_D^*$ denote the associated
dual eigenvalue sequence.  Comparing 
(\ref{DRECUR}) and (\ref{FIPOLY2}) using
(\ref{KI}) we routinely obtain
\begin{equation}\label{FITHETA1}
f_i(\theta)= k_i \theta^*_i/ \theta^*_0 \qquad \qquad (0 \leq i \leq 
D).
\end{equation}
We remark on two special cases. Setting $i=0,\; i=1$ in (\ref{FIPOLY2}) we 
routinely find $f_1=\lambda$, $\; f_{2}=(\lambda^{2}-k)/c_{2}$. Now setting 
$i=1, \; i=2$
in (\ref{FITHETA1}) we find
\begin{equation}\label{recth}
 \theta_1^*/\theta_0^* = \theta/k, \qquad \qquad 
 \theta_2^*/\theta_0^* = (\theta^{2}-k)/(kb_{1}).
\end{equation}

\bigskip
\noindent We now recall some polynomials related to the $f_i$. 
Let $p_0, p_1, \ldots, p_D$ denote the polynomials in ${\R}[\lambda]$ 
satisfying 
\begin{equation}
{p_{i} = \cases{f_{0}+f_{2}+f_{4}+\cdots +f_{i}, & if $i$ is even \cr
f_{1}+f_{3}+f_{5}+\cdots +f_{i}, & if $i$ is odd \cr}} 
 \qquad \qquad (0 \leq i \leq D).
\label{PIPOLY}
\end{equation}

\bigskip
\noindent
Let $i$ denote an integer $(0 \leq i \leq D)$.
The polynomial $p_i$ has degree $i$, and the
coefficient of $\lambda^i$ is
$(c_1c_2\cdots c_i)^{-1}$.  Recalling $f_{j}(A)=A_{j} \; (0 \le j \le 
D)$, we observe

\begin{equation}  \label{PDAeqs}
     p_{D}(A) + p_{D-1}(A) = J, \qquad \qquad p_{D}(A) - 
     p_{D-1}(A) = (-1)^{D}J',
\end{equation}
where $J'$ is from (\ref{J-def}).  

\bigskip
\noindent A bit later we find an orthogonality relation satisfied by
the polynomials $p_i$. To obtain it we use the following result.

\begin{lemma}\label{PIPROP}
Let $\Gamma$ denote a bipartite distance-regular graph with diameter $D 
\ge 4.$ 
Let the polynomials $f_0,f_1,\ldots,f_D$ be from (\ref{FIPOLY2}), and 
let
the polynomials $p_0,p_1,\ldots,p_D$ be from (\ref{PIPOLY}). Then
\begin{description} \item[(i)] $p_i - p_{i-2} = f_i \qquad \qquad (2
\le i \le D)$,
\item[(ii)] $(k^{2} - \lambda^{2})p_i = b_ib_{i+1}f_i - c_{i+1}c_{i+2}f_{i+2} \qquad \qquad 
(0 \le i \le D-2).$
\end{description} 
\end{lemma}
{\it Proof.} (i) Immediate from (\ref{PIPOLY}).\\
(ii) To see that the two sides are equal, in the expression on
the left
eliminate $p_i$ using
(\ref{PIPOLY}), and evaluate the result by repeatedly applying
(\ref{FIPOLY2}) and
Lemma \ref{bipequiv}(iv).
\hfill $\Box$\\

\begin{theorem}\label{PIRECUR}
Let $\Gamma$ denote a bipartite distance-regular graph with diameter $D 
\ge 4.$ 
Let the polynomials $p_0,p_1,\ldots,p_D$ be as in (\ref{PIPOLY}). 
Then $p_0=1$ and 
\begin{equation}  \label{PIRECUR1}
\lambda p_i = c_{i+1}p_{i+1}+b_{i+1}p_{i-1}
\qquad (0 \leq i \leq D-1),
\end{equation}
where $p_{-1}=0$.  
\end{theorem}
{\it Proof.}  Evaluate each 
side of (\ref{PIRECUR1}) using 
(\ref{PIPOLY}) and (\ref{FIPOLY2}), and simplify using Lemma 
\ref{bipequiv}(iv).
\hfill $\Box$\\

\noindent The polynomials $p_i$ satisfy the following orthogonality
relation.

\begin{lemma}\label{PIPJ}
Let $\Gamma=(X,R)$ denote a bipartite distance-regular graph with diameter $D 
\ge 4$ and eigenvalues $k=\theta_0 >\theta_1 >\cdots >\theta_D.$ Let 
the polynomials $p_0,p_1,\ldots, p_D$ be as in (\ref{PIPOLY}).
Then $p_D(\theta_h)=0$ and $p_{D-1}(\theta_h)=0$  for $1 \leq h \leq D-1$.
Moreover,
\begin{equation}
\label{PIPJ1}
\sum_{h=0}^D p_i(\theta_h)p_j(\theta_h)(k^{2}-\theta_h^{2})m_h = \delta_{ij} 
|X|k_ib_ib_{i+1}
\qquad \qquad (0 \le i,j \le D-2).
\end{equation}
(We recall $m_h$ denotes the multiplicity of $\theta_h$ for $0 \leq h 
\leq D$.)
\end{lemma}
{\it Proof.}  We first show $p_D(\theta_h)=0$ and 
$p_{D-1}(\theta_h)=0$  for $1 \leq h \leq D-1$.  Let $h$ be given.  
Multiplying both sides of the equations in (\ref{PDAeqs}) by $E_{h}$ 
and recalling $J, J'$ are scalar multiples of $E_{0}, E_{D}$, 
respectively, we find
\begin{equation}  \label{PIPJ1.5}
    p_{D}(A)E_{h} + p_{D-1}(A)E_{h} =0, \qquad \qquad 
    p_{D}(A)E_{h} - p_{D-1}(A)E_{h} =0
\end{equation}
  in view of (ev).  Solving the equations in (\ref{PIPJ1.5}), 
  we find $p_{D}(A)E_{h} =0,$ $\; p_{D-1}(A)E_{h} =0$.  
  Observe $p_{D}(A)E_{h} 
  =p_{D}(\theta_{h})E_{h}$, $\;p_{D-1}(A)E_{h} 
  =p_{D-1}(\theta_{h})E_{h}$, and thus $p_D(\theta_h)=0$, $p_{D-1}(\theta_h)=0$.
  Concerning (\ref{PIPJ1}), let the integers 
  $i,j$ be given. Without 
loss of generality,
we may assume $i \le j.$  We first assume $i$ is even.  
By (\ref{PIPOLY}) and Lemma \ref{PIPROP}(ii), 
the left-hand side of (\ref{PIPJ1}) is equal to
\begin{equation}\label{PIPJ2}
\sum_{h=0}^D \left(f_i(\theta_h)+ f_{i-2}(\theta_h)+\cdots 
+f_{0}(\theta_{h}) \right)
\left(b_j b_{j+1}f_j(\theta_h)- c_{j+1}c_{j+2}f_{j+2}(\theta_h)\right)m_h.
\end{equation}
Evaluating (\ref{PIPJ2}) using (\ref{FIFJ1}) and recalling $i \le j,$ 
we find
(\ref{PIPJ2}) is equal to the right-hand side of (\ref{PIPJ1}). The 
result follows.  The case in which $i$ is odd is similar.  \hfill $\Box$\\

\noindent The following fact will be useful.

\begin{lemma}\label{PITHETA*}
Let $\Gamma$ denote a bipartite distance-regular graph with diameter $D \ge
4$ and eigenvalues $k=\theta_0 > \theta_1 > \cdots > \theta_D$.  
Let $\theta$ denote one of $\theta_{1}, \theta_{2}, \ldots, 
\theta_{D-1}$ and let $\theta_0^*, \theta_{1}^{*},
\ldots, \theta_D^*$ denote
the corresponding dual eigenvalue sequence.  Then $\theta_0^*
\not= \theta_2^*$ and
\begin{equation}\label{PITHETA1}
p_i(\theta) = \displaystyle{\frac{b_2b_3\cdots b_{i+1}}{c_1c_2\cdots c_i}
\,\frac{\theta_i^*-\theta_{i+2}^*}{\theta_0^*-\theta_2^*}}
\qquad \qquad (0 \le i \le D-2),
\end{equation}
where the polynomials $p_i$ are from
(\ref{PIPOLY}).
\end{lemma}
{\it Proof.} Using the equation on the right in (\ref{recth}), we 
routinely verify that $\theta_0^*
\not= \theta_2^*$.  To obtain (\ref{PITHETA1}),
set $\lambda=\theta$ in Lemma \ref{PIPROP}(ii), and simplify the 
result using (\ref{KI}),
(\ref{FITHETA1}), and (\ref{recth}).
\hfill $\Box$\\

\begin{corollary}  \label{thfrac}
Let $\Gamma$ denote a bipartite distance-regular graph with diameter $D \ge
4$ and eigenvalues $k=\theta_0 > \theta_1 > \cdots > \theta_D$.  
Let $\theta$ denote one of $\theta_{1}, \theta_{2}, \ldots, 
\theta_{D-1}$ and let $\theta_0^*, \theta_{1}^{*},
\ldots, \theta_D^*$ denote
the corresponding dual eigenvalue sequence.
Then 
\begin{equation}  \label{sally}
    \frac{\theta^{2}-b_{2}}{b_{2}b_{3}} = \frac{\theta_2^*-\theta_4^*}
    {\theta_0^*-\theta_2^*}.
\end{equation}
We remark that the denominator on the right in (\ref{sally}) is 
nonzero by Lemma \ref{PITHETA*}.  
\end{corollary} 
{\it Proof.}   Using (\ref{PIRECUR1}) and recursion, we routinely find 
that $p_{2}(\theta) = (\theta^{2}-b_{2})/c_{2}$.  Using this fact and 
setting $i=2$ in (\ref{PITHETA1}), we obtain the desired result.
\hfill $\Box$\\


\section{The polynomials $\Psi_{i}$}

\bigskip
\noindent  Let $\Gamma$ denote a bipartite distance-regular graph with 
diameter $D \geq 4$.  In the previous section we used $\Gamma$ to 
define two families of polynomials in one variable.  We called these 
polynomials the $f_{i}$ and the $p_{i}$.  Later in this paper we will 
use $\Gamma$ to define a third family of polynomials in one 
variable.  We will call these polynomials the $g_{i}$.  To define and 
study the $g_{i}$ it is convenient to first consider some 
polynomials $\Psi_{i}$ in two variables.

\begin{definition} \label{PSIPOLY}  \rm
Let $\Gamma$ denote a bipartite distance-regular graph with diameter $D 
\ge 4$.  For $0 \leq i \leq D-2$ let $\Psi_{i}$ denote the polynomial 
in $\R[\lambda, \mu]$ given by
\begin{equation}
\Psi_i = \sum_{{h=0}\atop{ i-h {\tiny \mbox{ even}}}}^{i} 
p_{h}(\lambda) p_{h}(\mu)
\frac{k_{i}b_{i}b_{i+1}}
{k_{h}b_{h}b_{h+1}},
\label{PSIPOLY1}
\end{equation}
where the polynomials $p_0, p_{1}, \ldots, p_{D-2}$ are from
(\ref{PIPOLY}).  We observe $\Psi_{0}=1$ and $\Psi_{1} = \lambda 
\mu$.  
\end{definition}

\begin{lemma}
Let $\Gamma$ denote a bipartite distance-regular graph with diameter 
$D \geq 4$.  Let the polynomials $p_{i}, \Psi_{i}$ be as in 
(\ref{PIPOLY}), (\ref{PSIPOLY1}), respectively.  Then
\begin{equation}  \label{pipsi}
   p_{i}(\lambda) p_{i}(\mu) = \Psi_{i} 
   -\frac{b_{i}b_{i+1}}{c_{i}c_{i-1}}\Psi_{i-2} \qquad  \qquad (2 
   \leq i \leq D-2).
\end{equation}
\end{lemma}
{\it Proof.}  This is immediate from Definition \ref{PSIPOLY}.
\hfill $\Box$\\

\bigskip \noindent
The following equation is a variation of the Christoffel-Darboux 
Formula.

\begin{lemma}  \label{CDlemma}
Let $\Gamma$ denote a bipartite distance-regular graph with diameter 
$D \geq 4$.  Let the polynomials $p_{i}, \Psi_{i}$ be as in 
(\ref{PIPOLY}), (\ref{PSIPOLY1}), respectively.  Then for $1 \leq i 
\leq D-1$,
\begin{equation}  \label{CD}
  p_{i+1}(\lambda) p_{i-1}(\mu) - p_{i-1}(\lambda)p_{i+1}(\mu) = 
  c_{i}^{-1} c_{i+1}^{-1}(\lambda^{2}-\mu^{2})\Psi_{i-1}.
\end{equation}
\end{lemma}
{\it Proof.}  Repeatedly applying (\ref{PIRECUR1}), we find that for 
$0 \leq h \leq D-2$,
\begin{equation} \label{CD1}
   \lambda^{2} p_{h}(\lambda) = c_{h+1}c_{h+2}p_{h+2}(\lambda) + 
   (c_{h+1}b_{h+2}+b_{h+1}c_{h})p_{h}(\lambda) + 
   b_{h}b_{h+1}p_{h-2}(\lambda),
\end{equation}
 where we define $p_{-1}:=0, p_{-2} :=0$.  Similarly,
 \begin{equation} \label{CD2}
   \mu^{2} p_{h}(\mu) = c_{h+1}c_{h+2}p_{h+2}(\mu) + 
   (c_{h+1}b_{h+2}+b_{h+1}c_{h})p_{h}(\mu) + 
   b_{h}b_{h+1}p_{h-2}(\mu).
\end{equation}
  Subtracting $(k_{h}b_{h}b_{h+1})^{-1}p_{h}(\lambda)$ times 
  (\ref{CD2}) from $(k_{h}b_{h}b_{h+1})^{-1}p_{h}(\mu)$ times 
  (\ref{CD1}) and using (\ref{KI}), we find
\begin{equation}  \label{CD3}
  \frac{(\lambda^{2}-\mu^{2})p_{h}(\lambda)p_{h}(\mu) }
  {k_{h}b_{h}b_{h+1}} = 
  \frac{p_{h+2}(\lambda)p_{h}(\mu) - 
  p_{h}(\lambda)p_{h+2}(\mu)}{k_{h+2}} - 
  \frac{p_{h}(\lambda)p_{h-2}(\mu) - p_{h-2}(\lambda)p_{h}(\mu)}{k_{h}},
\end{equation}
   for $0 \leq h \leq D-2$.  Fix an 
   integer $i$ $(1 \leq i \leq D-1)$.  Summing  (\ref{CD3}) over all 
   $h$ such that $0 \leq h \leq i-1$ and such that $i-1-h$ is even, 
   and using (\ref{KI}), (\ref{PSIPOLY1}), we obtain (\ref{CD}).  \hfill $\Box$\\
   
\begin{lemma}  \label{PSIIPSIJ}
Let $\Gamma =(X,R)$ denote a bipartite distance-regular graph with diameter $D 
\ge 4$ and eigenvalues $k=\theta_0 >\theta_1 >\cdots >\theta_D.$  Let 
the polynomials $p_{i}, \; \Psi_{i}$ be as in (\ref{PIPOLY}), 
(\ref{PSIPOLY1}), respectively.  Then for $0 \leq i,j \leq D-2$,
\begin{equation} \label{PSIIPSIJ1}
\sum_{h=0}^{D}\Psi_{i}(\theta_{h},\mu)\Psi_{j}(\theta_{h},\mu)(k^{2}-\theta_{h}^{2})
(\mu^{2}-\theta_{h}^{2})m_{h} 
\end{equation}
\begin{equation} \label{PSIIPSIJ2}
  =\delta_{ij}|X|p_{i}(\mu)p_{i+2}(\mu)k_{i}b_{i}b_{i+1}c_{i+1}c_{i+2}.
 \end{equation}
 (We recall $m_{h}$ denotes the multiplicity of $\theta_{h}$ for $0 
 \leq h \leq D.$)
\end{lemma}
{\it Proof.}  Let the integers $i,j$ be given.  Without loss of 
generality we may assume $i \leq j$.  By Lemma \ref{CDlemma} we have
 $$ p_{j+2}(\theta_{h})p_{j}(\mu) - p_{j}(\theta_{h})p_{j+2}(\mu) = 
  c_{j+1}^{-1}c_{j+2}^{-1}(\theta_{h}^{2}-\mu^{2}) 
  \Psi_{j}(\theta_{h},\mu) $$
for $0 \leq h \leq D,$ so the sum (\ref{PSIIPSIJ1}) is equal to 
\begin{equation}  \label{psiorthog3}
c_{j+1}c_{j+2}\sum_{h=0}^{D}\Psi_{i}(\theta_{h},\mu)(p_{j+2}(\mu)p_{j}(\theta_{h})-p_{j}(\mu)p_{j+2}
(\theta_{h}))(k^{2}-\theta_{h}^{2})m_{h}.
\end{equation}
Eliminating $\Psi_{i}(\theta_{h}, \mu)$ in (\ref{psiorthog3}) using 
(\ref{PSIPOLY1}), and evaluating the result using Lemma \ref{PIPJ} 
and $i \leq j$, we obtain (\ref{PSIIPSIJ2}).  The result follows.  \hfill $\Box$\\

\begin{lemma}  \label{PITHETA}
Let $\Gamma$ denote a bipartite distance-regular graph with diameter $D 
\ge 4$ and eigenvalues $k=\theta_0 >\theta_1 >\cdots >\theta_D.$  Let 
the polynomials $p_{i}$ be as in (\ref{PIPOLY}). 
Then the following (i), (ii) hold. 
\begin{description}
\item[(i)]  Abbreviate $\theta=\theta_{1}$.  Then $p_{i}(\theta) > 0 $
for $0 \le i \le D-2$, and $p_{D-1}(\theta)=0, \, p_{D}(\theta)=0$.
\item[(ii)]  Abbreviate $\theta=\theta_{D-1}$.  Then $(-1)^{i}p_{i}(\theta) > 
0 $
for $0 \le i \le D-2$, and $p_{D-1}(\theta)=0, \, p_{D}(\theta)=0$.
\end{description}
\end{lemma}
{\it Proof.}  (i)  Observe $p_{D-1}(\theta_{1})=0, \, 
p_{D}(\theta_{1})=0$ by Lemma \ref{PIPJ}.
Suppose there exists an integer $i$ $(0 \le i \le D-2)$ 
such that $p_{i}(\theta) \le 0$.  Let us pick the minimal such $i$.  
Observe $i \ge 2$ since $p_{0}(\theta) =1$, $p_{1}(\theta) = 
\theta$.  Apparently $p_{i-2}(\theta) > 0$.  We claim there exists an 
integer $h$ $(2 \le h \le D-2)$ such that $\Psi_{i-2}(\theta_{h}, 
\theta) \not= 0$.  To see this, observe by Definition \ref{PSIPOLY} 
that $\Psi_{i-2}(\lambda, \theta)$ is a polynomial in $\lambda$ with 
degree $i-2$.  In this polynomial the coefficient of $\lambda^{i-2}$ 
is $p_{i-2}(\theta)(c_{1}c_{2}\cdots c_{i-2})^{-1}$.  Apparently this 
polynomial is not identically 0 so there exist at most $i-2$ integers 
$h$ $(2 \le h \le D-2)$ such that $\Psi_{i-2}(\theta_{h}, 
\theta)=0$.  By this and since $i \le D-2$, there exists at least one 
integer $h$ $(2 \le h \le D-2)$ such that $\Psi_{i-2}(\theta_{h}, 
\theta) \not= 0$.  We have now proved our claim.  We may now argue
\begin{eqnarray*}
 0 &<& \sum_{h=2}^{D-2} 
 \Psi_{i-2}^{2}(\theta_{h},\theta)(k^{2}-\theta_{h}^{2}) 
 (\theta^{2}-\theta_{h}^{2})m_{h} \\
 &=& \sum_{h=0}^{D} 
 \Psi_{i-2}^{2}(\theta_{h},\theta)(k^{2}-\theta_{h}^{2}) 
 (\theta^{2}-\theta_{h}^{2})m_{h}  \\
&=& |X|p_{i-2}(\theta)p_{i}(\theta)k_{i-2}b_{i-2}b_{i-1}c_{i-1}c_{i} 
 \quad
 \mbox{ (by Lemma \ref{PSIIPSIJ}) }\\ & \le& 0.
 \end{eqnarray*}
 We now have a contradiction and the result follows. \\
 (ii) Similar to the proof of (i) above.  \hfill $\Box$\\

\begin{lemma}  \label{PITHETAD}
Let $\Gamma$ denote a bipartite distance-regular graph with diameter $D 
\ge 4$ and eigenvalues $k=\theta_0 >\theta_1 >\cdots >\theta_D.$  
Assume $D$ is odd and abbreviate $d = (D-1)/2$.  Let 
the polynomials $p_{i}$ be as in (\ref{PIPOLY}). 
Then the following (i), (ii) hold.
\begin{description}
\item[(i)]  Abbreviate $\theta=\theta_{d}$.  Then $(-1)^{\lfloor 
\frac{i}{2}\rfloor}p_{i}(\theta) > 0 $
for $0 \le i \le D-2$, and $p_{D-1}(\theta)=0, \, p_{D}(\theta)=0$.
\item[(ii)]  Abbreviate $\theta=\theta_{d+1}$.  Then $(-1)^{\lfloor 
\frac{i+1}{2}\rfloor}p_{i}(\theta) > 0 $
for $0 \le i \le D-2$, and $p_{D-1}(\theta)=0, \, p_{D}(\theta)=0$.
\end{description}
\end{lemma}
{\it Proof.}  (i)  Observe $p_{D-1}(\theta_{d})=0, \, 
p_{D}(\theta_{d})=0$ by Lemma \ref{PIPJ}.  Suppose there exists an integer 
$i$ $(0 \le i \le D-2)$ 
such that $(-1)^{\lfloor \frac{i}{2}\rfloor}p_{i}(\theta) \le 0$.  
Let us pick the minimal such $i$.  
Observe $i \ge 2$ since $p_{0}(\theta) =1$, $p_{1}(\theta) = 
\theta$.  Apparently $(-1)^{\lfloor \frac{i-2}{2}\rfloor}p_{i-2}(\theta) > 
0$, so $p_{i}(\theta), \, p_{i-2}(\theta)$ do not have opposite 
signs.    
We claim there exists an 
integer $h$ $(1 \le h \le D-1), \; h\not= d,\, h \not= d+1$ such that 
$\Psi_{i-2}(\theta_{h}, 
\theta) \not= 0$.  To see this, observe by Definition \ref{PSIPOLY} 
that $\Psi_{i-2}(\lambda, \theta)$ is a polynomial in $\lambda$ with 
degree $i-2$.  In this polynomial the coefficient of $\lambda^{i-2}$ 
is $p_{i-2}(\theta)(c_{1}c_{2}\cdots c_{i-2})^{-1}$.  Apparently this 
polynomial is not identically 0 so there exist at most $i-2$ integers 
$h$ $(1 \le h \le D-1)$ such that $\Psi_{i-2}(\theta_{h}, 
\theta)=0$.  By this and since $i \le D-2$, there exists at least one 
integer $h$ $(1 \le h \le D-1), \; h\not= d,\, h \not= d+1$ 
such that $\Psi_{i-2}(\theta_{h}, 
\theta) \not= 0$.  We have now proved our claim.  We may now argue
\begin{eqnarray*}
 0 &>& \sum_{{h=1}\atop{h\not= d, h\not= d+1}}^{D-1} 
 \Psi_{i-2}^{2}(\theta_{h},\theta)(k^{2}-\theta_{h}^{2}) 
 (\theta^{2}-\theta_{h}^{2})m_{h} \\
 &=& \sum_{h=0}^{D} 
 \Psi_{i-2}^{2}(\theta_{h},\theta)(k^{2}-\theta_{h}^{2}) 
 (\theta^{2}-\theta_{h}^{2})m_{h}  \\
 &=& |X|p_{i-2}(\theta)p_{i}(\theta)k_{i-2}b_{i-2}b_{i-1}c_{i-1}c_{i} 
 \quad
 \mbox{ (by Lemma \ref{PSIIPSIJ}) }\\ & \ge& 0.
 \end{eqnarray*}
 We now have a contradiction and the result follows. \\
 (ii)  Similar to the proof of (i) above.  \hfill $\Box$\\


\section{A third family of polynomials}

\bigskip
\noindent In this section we will use the following notation.


\begin{notation} \label{thetanote}  \rm
Throughout this section, $\Gamma$ will denote a bipartite distance-regular graph with diameter $D 
\ge 4$ and eigenvalues $k=\theta_0 >\theta_1 >\cdots >\theta_D.$  Let 
the polynomials $p_{i}$ be as in (\ref{PIPOLY}). 
If $D$ is odd, we let $\theta$ 
denote
one of $\theta_1, \theta_d, \theta_{d+1}, \theta_{D-1}$,
where $d=(D-1)/2$.    
If $D$ is even, we let $\theta$ denote
one of $\theta_{1}, \theta_{D-1}$.  We remark 
$p_{i}(\theta) \ne 0$ for $0 \le i \le D-2$ by Lemma \ref{PITHETA} 
and Lemma \ref{PITHETAD}.
\end{notation}

\bigskip
\noindent In
this section we use $\Gamma$
to define a family of polynomials in one variable. We call these 
polynomials
the $g_i$.

\begin{definition}\label{GIPOLY}  \rm
With reference to Notation \ref{thetanote},
for $0 \leq i \leq D-2$ we define the polynomial
$g_i \in \R[\lambda]$ by
$$g_i = \sum_{{h=0}\atop{ i-h {\tiny \mbox{ even}}}}^{i}
\frac{p_h(\theta)k_i b_i b_{i+1}}{p_i(\theta) k_hb_h b_{h+1}} p_h.$$
We observe 
\begin{equation} \label{PSIG}
\Psi_{i}(\lambda, \theta) = p_{i}(\theta)g_{i}, 
\end{equation}
where $\Psi_{i}$ is from Definition \ref{PSIPOLY}.  
We emphasize $g_i$ depends on $\theta$ as well as
the intersection numbers of $\Gamma$.
\end{definition}

\begin{lemma}  
With reference to Notation \ref{thetanote}, let 
$g_{0}, g_{1}, \ldots, g_{D-2}$ denote the associated polynomials from 
Definition \ref{GIPOLY}.  Then
\begin{equation}  \label{pigi}
  p_{i} = g_{i} - 
  \frac{b_{i}b_{i+1}}{c_{i-1}c_{i}}\frac{p_{i-2}(\theta)}{p_{i}(\theta)}g_{i-2} 
   \qquad (2 \le i \le D-2).
\end{equation} 
\end{lemma}
{\it Proof.}  Set $\mu = \theta$ in (\ref{pipsi}) and simplify the 
result using (\ref{PSIG}).  \hfill $\Box$\\

\begin{lemma}  
With reference to Notation \ref{thetanote}, let 
$g_{0}, g_{1}, \ldots, g_{D-2}$ denote the associated polynomials from 
Definition \ref{GIPOLY}.  Then (i) and (ii) hold below for $0 \le i 
\le D-2$:
\begin{description}
\item[(i)]  The polynomial $g_{i}$ has degree exactly $i$.  
\item[(ii)]  The coefficient of $\lambda^{i}$ in $g_{i}$ is 
$(c_{1}c_{2}\cdots c_{i})^{-1}$.  
\end{description}
\end{lemma}
{\it Proof.}  Routine.  \hfill $\Box$\\

\bigskip \noindent
We now present a three-term recurrence satisfied by the polynomials $g_{i}$.

\begin{theorem}  \label{Grecurthm}  
With reference to Notation \ref{thetanote}, let 
$g_{0}, g_{1}, \ldots, g_{D-2}$ denote the associated polynomials from 
Definition \ref{GIPOLY}.  Then $g_{0}=1$ and
\begin{equation}  \label{girecur}
  \lambda g_{i} = c_{i+1}g_{i+1} + \omega_{i} g_{i-1} 
\end{equation}
   for $0 \le i \le D-3$, where $g_{-1}=0, \omega_{0}=0$,
and
\begin{equation}  \label{omega}
  \omega_{i} = \frac{b_{i+1}c_{i+2}}{c_{i}} \frac{p_{i-1}(\theta) 
  p_{i+2}(\theta)}{p_{i}(\theta) p_{i+1}(\theta)} \qquad \qquad (1 
  \le i \le D-3).
\end{equation}
\end{theorem}
{\it Proof.} We find $g_{0}=1$ by Definition \ref{GIPOLY}.  To 
prove (\ref{girecur}), we proceed by induction.  
It is routine to show equality holds in (\ref{girecur}) for 
$i=0,1$ using Definition \ref{GIPOLY} and (\ref{PIRECUR1}).  
Now suppose equality holds in (\ref{girecur}) for $i=j-2$, where $ 2 
\le j  \le D-3$.  We show equality holds in (\ref{girecur}) 
for $i=j$.  By the inductive hypothesis, we have
\begin{equation}  \label{girecurA}
  \lambda g_{j-2} =  c_{j-1} g_{j-1} + \omega_{j-2} g_{j-3} .
\end{equation}
We must prove $\lambda g_{j} = c_{j+1}g_{j+1}+\omega_{j}g_{j-1}$.  
First consider the expression $c_{j+1}g_{j+1} + \omega_{j} g_{j-1}$.
Eliminating $g_{j+1}, \omega_{j}$ in this expression 
using (\ref{pigi}), (\ref{omega}), respectively, and then simplifying 
the result using (\ref{PIRECUR1}), we 
find  
\begin{equation}   \label{girecurD}
  c_{j+1}g_{j+1} + \omega_{j} g_{j-1}
   = c_{j+1} p_{j+1} + \frac{b_{j+1} \theta p_{j-1}(\theta)}{c_{j} 
   p_{j}(\theta)} g_{j-1}.
\end{equation}
Now consider the expression $\lambda g_{j}$.  Replacing $g_{j}$ in 
this expression  using (\ref{pigi}), and eliminating $\lambda p_{j}$, 
$\lambda g_{j-2}$ in the result using (\ref{PIRECUR1}), 
(\ref{girecurA}), respectively, we find
\begin{equation}  \label{girecurB}
   \lambda g_{j} =  c_{j+1}p_{j+1} + b_{j+1} p_{j-1} + 
   \frac{b_{j}b_{j+1}}{c_{j-1}c_{j}} 
   \frac{p_{j-2}(\theta)}{p_{j}(\theta)} (c_{j-1} g_{j-1} + 
   \omega_{j-2} g_{j-3}). 
\end{equation}
  Eliminating $\omega_{j-2}$ in (\ref{girecurB}) using (\ref{omega}) 
  and eliminating 
  $b_{j-1}b_{j}p_{j-3}(\theta)(c_{j-2}c_{j-1}p_{j-1}(\theta))^{-1} 
  g_{j-3}$ in the result using (\ref{pigi}), we find
\begin{equation}  \label{girecurC}
 \lambda g_{j} = c_{j+1}p_{j+1} + b_{j+1} 
 \, \frac{c_{j}p_{j}(\theta) 
 +b_{j}p_{j-2}(\theta)}{c_{j}p_{j}(\theta)} \, g_{j-1}.
\end{equation}
Observe the right-hand sides of (\ref{girecurD}), (\ref{girecurC}) are 
equal
in view of (\ref{PIRECUR1}), and thus the left-hand sides are equal.
We find equality holds in (\ref{girecur}) 
for $i=j$, as desired.  \hfill $\Box$\\

\begin{theorem} \label{GIGJ}
With reference to Notation \ref{thetanote}, let 
$g_{0}, g_{1}, \ldots, g_{D-2}$ denote the associated polynomials from 
Definition \ref{GIPOLY}.  Then
\begin{equation}  \label{GIGJ1}
\sum_{h=0}^{D}g_{i}(\theta_{h})g_{j}(\theta_{h})(k^{2}-\theta_{h}^{2})
(\theta^{2}-\theta_{h}^{2}) m_{h} = 
\delta_{ij}|X|k_{i}b_{i}b_{i+1}c_{i+1}c_{i+2}\frac{p_{i+2}(\theta)}{p_{i}(\theta)}
\end{equation}
for $0 \le i,j \le D-2$.  
\end{theorem}
{\it Proof.}  To verify (\ref{GIGJ1}), on the left-hand side first 
eliminate $g_{i}(\theta_{h})$ and $g_{j}(\theta_{h})$ using 
(\ref{PSIG}), and then evaluate the result using Lemma 
\ref{PSIIPSIJ}.  \hfill $\Box$\\


\section{The subconstituent algebra and its modules}
In this section we recall some definitions and basic concepts 
concerning
the subconstituent algebra and its modules.
For more information we refer the reader to \cite{caugh2},
\cite{curtin1},
\cite{curtin2},
\cite{go},
\cite{hobart},
\cite{terwSub1}.

\bigskip
\noindent
Let $\Gamma=(X,R)$ denote a distance-regular graph with diameter $D 
\ge 3.$
We recall the dual Bose-Mesner algebra of $\Gamma.$
For the rest of this section, fix a vertex $x \in X.$ For $ 0 \le i 
\le D$ let $E_i^*=E_i^*(x)$ denote the diagonal
matrix in $\hbox{Mat}_X(\C)$ with $yy$ entry
\begin{equation}\label{DEFDEI}
{(E_i^*)_{yy} = \cases{1, & if $\partial(x,y)=i$\cr
0, & if $\partial(x,y) \ne i$\cr}} \qquad (y \in X).
\end{equation}
We call $E_i^*$ the {\it $i^{\hbox {th}}$ dual idempotent of} $\Gamma$
{\it with respect to x.} 
We observe
(di) $\sum_{i=0}^D E_i^* = I$;
(dii) $\overline{E_i^*} = E_i^* \; (0 \le i \le D)$; 
(diii) $E_i^{*t} = E_i^*  \;(0 \le i \le D)$;
(div) $E_i^*E_j^* = \delta_{ij}E_i^*  \;(0 \le i,j \le 
D).$
Using (di) and (div) 
we find $E_0^*,E_1^*, \ldots, E_D^*$ form a basis for a commutative 
subalgebra $M^*=M^*(x)$ of $\hbox{Mat}_X(\C).$ We call $M^*$ the {\it 
dual Bose-Mesner algebra of}
$\Gamma$ {\it with respect to x.} 
We recall the subconstituents of $\Gamma $.
Using (\ref{DEFDEI}) we find
\begin{equation}\label{DEIV}
E_i^*V = {\rm Span}\,\{\hat{y} \;|\; y \in X, \quad \partial(x,y)=i\}
\qquad (0 \le i \le D).
\end{equation}
By (\ref{DEIV}) and since $\lbrace {\hat y} \;|\;y \in X\rbrace $ is
an orthonormal basis for $V$ we find 
\begin{equation}
\label{EISORTHO}
V = E_0^*V+E_1^*V+ \cdots +E_D^*V \qquad \qquad {\rm (orthogonal\ 
direct\ sum}).
\end{equation}
Combining (\ref{DEIV}) and (\ref{DEFKI}) we find
the dimension of $E_i^*V$ is $k_i$ for $(0 \le i \le D)$.
We call $E_i^*V$ the {\it $i^{\hbox {th}}$ subconstituent of} $\Gamma$
{\it with respect to} $x$.

\bigskip
\noindent
We recall how $M$ and $M^*$ are related.
By \cite[Lemma 3.2]{terwSub1}, \begin{equation} \label{REL2}
E_h^*A_iE_j^*=0 \quad {\rm if\ and\ only\ if} \quad p_{ij}^h = 0
\qquad \qquad (0 \le h,i,j \le D).
\end{equation}

\noindent
Let $T=T(x)$ denote the subalgebra of $\hbox{Mat}_X(\C)$ generated by 
$M$ and $M^*$. We call $T$ the {\it subconstituent algebra of} 
$\Gamma$ {\it with respect to} $x$ \cite{terwSub1}. We observe
$T$ has finite dimension. Moreover $T$ is semi-simple;
the reason is that $T$ is closed under the conjugate-transpose map
\cite[p. 157]{CR}.

\bigskip
\noindent
We now consider the modules for $T.$ By a {\it T-module}
we mean a subspace $W \subseteq V$ such that $BW \subseteq W$
for all $B \in T.$ We refer to $V$ itself as the {\it standard
module} for $T.$ Let $W$ denote a $T$-module. Then $W$ is said
to be {\it irreducible} whenever $W$ is nonzero and $W$ contains no 
$T$-modules other than 0 and $W.$ Let $W,W^\prime$ denote
$T$-modules. By an {\it isomorphism of $T$-modules}
from $W$ to $W^\prime$ we
mean an isomorphism of vector spaces
$\sigma: W \rightarrow W^\prime$
such that
$$
(\sigma B- B \sigma)W = 0 \qquad \qquad {\rm for\ all}\; B \in T.
$$
The modules $W,W^\prime$ are said to be {\it isomorphic as 
$T$-modules}
whenever
there exists an isomorphism of $T$-modules from $W$ to $W^\prime.$

\bigskip
\noindent
Let $W$ denote a $T$-module and let $W'$ denote a $T$-module 
contained in $W$. Using (\ref{ADJ}) we find the orthogonal complement 
of $W'$ in $W$ is a $T$-module.
It follows that each $T$-module
is an orthogonal direct sum of irreducible $T$-modules.
We mention any two nonisomorphic
irreducible $T$-modules are orthogonal \cite[Chapter IV]{CR}.

\bigskip
\noindent
Let $W$ denote an irreducible $T$-module.
Using (di)--(div) above we find $W$ is the direct sum of the
nonzero spaces among $E^*_0W, E^*_1W,\ldots, E^*_DW$.
Similarly using
(eii)--(ev) we find
$W$ is the direct sum of the
nonzero spaces among $E_0W, E_1W,\ldots, E_DW$.
If the dimension of $E^*_iW$ is at most 1 for
$0 \leq i \leq D$ then
the dimension of $E_iW$ is at most 1 for
$0 \leq i \leq D$ \cite[Lemma 3.9]{terwSub1}; in this case we say $W$ 
is {\it thin}.
Let $W$ denote an irreducible $T$-module.
By the {\it endpoint}
of $W$ we mean $$ {\rm min}\,\{i\;|\;0\le i \le D, \;\; E_i^*W \ne 0 
\}.$$

\noindent
In the rest of the paper we will assume $\Gamma$ is bipartite.  
We adopt the following
notational convention.

\begin{definition}\label{A}  \rm
For the rest of this paper we let $\Gamma=(X,R)$ denote a bipartite
distance-regular graph with diameter $D \ge 4$, valency $k \ge 3$, 
intersection numbers $b_i, c_i$,
adjacency matrix $A$,
Bose-Mesner algebra $M$,
and eigenvalues
$\theta_0 > \theta_1 > \cdots > \theta_D$.
For $0 \leq i \leq D$ we let $E_i $ denote the primitive idempotent of
$\Gamma$ associated with $\theta_i$. We define $d = \lfloor D/2 \rfloor$.
We let $V$ denote the standard 
module for $\Gamma$. We fix $x \in X$ and abbreviate
$E_i^*=E_i^*(x)$ $(0 \le i \le D)$,
$ M^*=M^*(x)$,
$T=T(x)$.  We define
\begin{equation}
s_i=\sum_{{y \in X}\atop {\partial(x,y)=i}} {\hat y}
\qquad \qquad (0 \leq i \leq D).
\label{si}
\end{equation}
\end{definition}


\section{The $T$-module of endpoint 0}

With reference to Definition \ref{A}, there exists a unique
irreducible $T$-module with endpoint 0 \cite[Proposition 8.4]{egge1}. We 
call this module
$V_0$. The module $V_0$ is described in \cite{curtin1}, \cite{egge1}.
We summarize some details below in order to 
motivate the results that
follow.

\bigskip
\noindent The module $V_0$ is thin. In fact
each of $E_iV_0$, $E^*_iV_0$
has dimension 1 for $0 \leq i \leq D$. We give two bases for $V_0$. 
The vectors 
$E_0{\hat x}, E_1{\hat x}, \ldots , E_D{\hat x}$
form a basis for $V_0$. These vectors are mutually orthogonal and 
$\|E_i {\hat x}\|^2 = m_i|X|^{-1}$ for $0 \leq i \leq D$.
To motivate the second basis we make some comments.  
For $0 \leq i \leq D$ we have
$s_i = A_i{\hat x}$.
Moreover $s_i=E^*_i\delta $, where $\delta=\sum_{y \in X} {\hat y}$. 
The vectors 
$s_0, s_1, \ldots, s_D$
form a basis for $V_0$.
These vectors are mutually orthogonal and
$\|s_i\|^2 = k_i$ for $0 \leq i \leq D$.
With respect to the basis $s_0, s_1, \ldots, s_D$ the matrix 
representing $A$ is
\begin{eqnarray*}
\left(\begin{array}{cccccc}
0 & b_0 & & & & {\bf 0}\\
c_1 & 0 & b_1 & & & \\
& c_2 & \cdot & \cdot & & \\
& & \cdot & \cdot & \cdot& \\
& & & \cdot & \cdot & b_{D-1} \\
{\bf 0} & & & & c_{D} & 0
\end{array} \right).
\end{eqnarray*}
The two bases for $V_{0}$ given above
are related as follows.
For $0 \leq i \leq D$ we have
$$s_i = \sum_{h=0}^{D}f_i(\theta_h) E_h{\hat x},$$
where the polynomial $f_i$ is from
(\ref{FIPOLY2}).


\section{$T$-modules of endpoint 1}

With reference to Definition \ref{A}, there exists, up to isomorphism, a unique
irreducible $T$-module with endpoint 1 \cite[Corollary 7.7]{curtin1}. 
We call this module $V_1$.
The module $V_1$ is described in \cite{curtin1}, \cite{gotight}.
We summarize some details below in order to 
motivate the results that
follow.

\bigskip \noindent
The module $V_{1}$ is thin with dimension $D-1$.
We give two bases for $V_1$. Let $v$ denote a nonzero
vector in $E^*_1V_1$. Then the sequence
$E_1v, E_{2}v, \ldots, E_{D-1}v $  
is a basis for $V_1$. These vectors are mutually
orthogonal and
$$||E_iv||^2 =  \frac{m_i(k^2-\theta^2_i)}{|X|k(k-1)} ||v||^2 \qquad 
\qquad (1 \leq i \leq D-1).$$
To motivate the second basis we make some
comments.
We have
$E^*_{i+1}A_{i}v=p_i(A)v$ for   $0 \leq i \leq D-2$
and $E^*_DA_{D-1}v=0$, where the polynomial $p_{i}$ is from (\ref{PIPOLY}).
The vectors
$E^*_{1}A_{0}v, E^{*}_{2}A_{1}v, \ldots, E^{*}_{D-1}A_{D-2}v$  
form a basis for $V_1$.
These vectors are mutually orthogonal
and $$||E^*_{i+1}A_{i}v||^2 =
\frac{b_2 \cdots b_{i+1}}{c_1 \cdots c_{i}} ||v||^2
 \qquad \qquad (0 \leq i \leq D-2).$$
With respect to the basis 
$E^*_{1}A_{0}v, E^{*}_{2}A_{1}v, \ldots, E^{*}_{D-1}A_{D-2}v$, the
matrix representing $A$ is
\begin{eqnarray*}
\left(\begin{array}{cccccc}
0 & b_2 & & & & {\bf 0}\\
c_1 & 0 & b_3 & & & \\
& c_2 & \cdot & \cdot & & \\
& & \cdot & \cdot & \cdot& \\
& & & \cdot & \cdot & b_{D-1} \\
{\bf 0} & & & & c_{D-2} & 0
\end{array} \right).
\end{eqnarray*}
The two bases for $V_{1}$ given above
are related as follows.
For $0 \leq i \leq D-2$ we have
$$E^*_{i+1}A_{i}v = \sum_{h=1}^{D-1} p_i(\theta_h)E_hv.$$

\bigskip \noindent
We comment that $V_1$ appears in $V$ with multiplicity
$k-1$.
We will need the following result.

\begin{corollary} \label{V1eigval} With reference to Definition \ref{A},
let $W$ denote an irreducible $T$-module with endpoint
$1$.  Observe $E^*_2W$ is an eigenspace for $E^*_2A_2E^*_2$.
The corresponding eigenvalue is $b_3-1$.
\end{corollary}
{\it Proof.}  The desired eigenvalue is the entry in the second row 
and second column of the matrix
representing $A_2$ with respect to the basis $E^*_{1}A_{0}v, E^{*}_{2}A_{1}v, \ldots, E^{*}_{D-1}A_{D-2}v$.  
To compute this entry, first set $i=1$ 
in (\ref{AAI1}) and observe that $c_2 A_2=A^2-kI$.  Using this fact 
and the above matrix display of $A$, we verify the specified matrix 
entry is $b_{3}-1$.
\hfill $\Box$\\


\section{The local eigenvalues}

Later in the paper we will consider the thin irreducible $T$-modules 
with endpoint 2.  In order to discuss these we introduce some 
parameters we call the local eigenvalues.  

\begin{definition}\label{SUBGRAPH}  \rm
With reference to Definition \ref{A},
we let $\Gamma_{2}^{2} = \Gamma_{2}^{2}(x) $ denote the
graph
$(\breve{X},\breve{R}),$ where
\begin{eqnarray}
\breve{X} &=& \{y \in X\;|\;\partial(x,y)=2\},\nonumber\\
\breve{R} &=& \{yz\;|\;y,z \in \breve{X},\, \partial(y,z)=2 \},\nonumber
\end{eqnarray}
where we recall $\partial$ denotes the path-length distance function 
for $\Gamma$.
The graph $\Gamma_{2}^{2}$ has 
exactly $k_{2}$ vertices, where $k_{2}$ is the second valency of
$\Gamma.$ Also, $\Gamma_{2}^{2}$ is regular with valency $p^{2}_{22}$.
We let $\breve {A}$ denote the adjacency matrix of $\Gamma_{2}^{2} $. The 
matrix $\breve {A}$ is symmetric with real entries; therefore $\breve 
{A}$ is diagonalizable with all eigenvalues real.
We let
$\eta_1, \eta_2 , \ldots , \eta_{k_{2}}$
denote the eigenvalues of $\breve {A}$.  We call 
$\eta_1, \eta_2 , \ldots , \eta_{k_{2}}$ the {\em local eigenvalues
of $\Gamma$ with respect to $x$}.  
\end{definition}

\noindent With reference to Definition \ref{A}, we consider the
second subconstituent $E^*_2V$. We recall 
the dimension of $E^*_2V$ is $k_{2}$.
Observe $E^*_2V$ is invariant under the action of $E^*_2A_{2}E^*_2$.
To illuminate this action we make an observation.
For an appropriate ordering of the vertices of $\Gamma$ we have $$
E^*_2A_{2}E^*_2 =
\left(\begin{array}{cc} \breve{A} & 0 \\ 0 & 0 \end{array} \right),
$$
where $\breve {A}$ is from
Definition \ref{SUBGRAPH}. Apparently the action of
$E^*_2A_{2}E^*_2$ on $E^*_2V$ is essentially the adjacency map
for $\Gamma_{2}^{2} $.
In particular the action of $E^*_2A_{2}¥E^*_2$ on $E^*_2V$ is 
diagonalizable
with eigenvalues $\eta_1, \eta_2, \ldots, \eta_{k_{2}}$.  We observe 
the vector $s_{2}¥$ from (\ref{si}) is contained in $E_{2}^{*}V$.  
One may easily show that $s_{2}$ is an eigenvector for $E^*_2A_{2}E^*_2$ 
with eigenvalue $p_{22}^{2}$.  
Let $v$ denote a vector in $E_{2}^{*}V$.
We observe the following are equivalent:  (i) $v$ is orthogonal to
$s_2$; (ii) $E_0v=0$; (iii) $Jv=0$; (iv) $E_{D}v=0$; (v) $J'v=0$.  
Let $V_{1}$ 
denote an irreducible $T$-module of endpoint $1$, and let $v$ 
denote a vector in $E_{2}^{*}V_{1}$.  By Corollary \ref{V1eigval}, 
$v$ is an eigenvector for $E^*_2A_{2}¥E^*_2$ with eigenvalue 
$b_{3}-1$.  
Reordering 
the eigenvalues if necessary, we have $\eta_{1} = p_{22}^{2}$ 
and $\eta_{i}= b_{3}-1$ $(2 \leq i \leq k)$.  For the rest of this 
paper we assume the local eigenvalues of $\Gamma$ are ordered in this 
way.

\bigskip \noindent
We now need some notation.    

\begin{definition} \label{U=V0+Y}  \rm With reference to Definition \ref{A}, 
let $Y$ denote the subspace of $V$ spanned by the
   irreducible $T$-modules with endpoint 1.  We define 
the set $U$ to be the orthogonal 
   complement of $E^*_2V_{0} + E^*_2Y$ in $E^{*}_{2}V$.
 \end{definition}
 

\begin{definition}  \label{phidef}  \rm
With reference to Definition \ref{A}, let $\Phi$ denote the set of 
distinct scalars among $\eta_{k+1}$,$\eta_{k+2}$,$\ldots$, 
$\eta_{k_{2}}$, where the $\eta_{i}$ are from Definition 
\ref{SUBGRAPH}.  For $\eta \in \R$ we let $\mbox{mult}_{\eta}$ denote the 
number of times $\eta$ appears among 
$\eta_{k+1}, \eta_{k+2}, \ldots, \eta_{k_{2}}$.  We observe 
$\mbox{mult}_{\eta} \not= 0$ if and only if $\eta \in \Phi$.  
\end{definition}

\noindent
Using (\ref{ADJ}) we find $U$ is invariant under $E^*_2A_{2}E^*_2$. 
Apparently the restriction of $E^*_2A_{2}E^*_2$ to $U$ is diagonalizable 
with eigenvalues 
$\eta_{k+1}, \eta_{k+2}, \ldots, \eta_{k_{2}}$.
For $\eta \in \R$ let $U_\eta $ denote the set consisting of those
vectors in $U$ that are eigenvectors for $E^*_2A_{2}¥E^*_2$ with 
eigenvalue $\eta $.
We observe $U_\eta $ is a subspace of $U$ with dimension 
$\mbox{mult}_{\eta}$.
We emphasize the following are equivalent:  (i) $\mbox{mult}_{\eta} 
\not= 0$; (ii) $U_{\eta} \not= 0$; (iii) $\eta \in \Phi$.  
By (\ref{ADJ}) and since $E^*_2A_{2}¥E^*_2$ is symmetric with real entries 
we find
\begin{equation}
U = \sum_{\eta \in \Phi} U_\eta \qquad \qquad (\hbox{orthogonal 
direct sum}).
\label{Ubreakdown}
\end{equation}

\bigskip
\noindent The following result will be useful.

\begin{lemma}\label{DEEDE}
With reference to Definition \ref{A}, let $E$ denote a primitive 
idempotent of $\Gamma$ and let 
$\theta_0^*,\theta_1^*,\ldots,\theta_D^*$ denote the corresponding 
dual eigenvalue sequence.
Then
\begin{equation}\label{DEEDE1}
|X| E_2^*EE_2^* = (\theta_0^*-\theta_4^*)E_2^* + 
(\theta_2^*-\theta_4^*)E^*_2A_{2}¥E^*_2 + \theta_4^*E^*_2JE^*_2.
\end{equation}
\end{lemma}
{\it Proof.} By (aii),
\begin{eqnarray}
E_2^*JE_2^*&=& E_2^*\left(\sum_{i=0}^D A_i \right)E_2^*\nonumber\\
{} &=& E_2^* + E_2^*A_2E_2^* + E_2^*A_4E_2^* \label{EQUA1}
\end{eqnarray}
in view of
(ai), 
(\ref{REL2}).
Using (\ref{PRIMIDEM}) we similarly find
\begin{eqnarray}
|X|E_2^*EE_2^* &=& E_2^*\left(\sum_{i=0}^D \theta_i^* A_i 
\right)E_2^*\nonumber\\
{} &=& \theta_0^*E_2^* + \theta_2^*E_2^*A_2E_2^* + 
\theta_4^*E_2^*A_4E_2^*. \label{EQUA2}
\end{eqnarray}
Eliminating $E_2^*A_4E_2^*$ in (\ref{EQUA2})
using (\ref{EQUA1})
we get (\ref{DEEDE1}).
\hfill $\Box$\\


\section{Some inner products}

\noindent
In this section we are concerned with the following basis.

\begin{lemma}\label{EIMVBASIS}
With reference to Definition \ref{A},
let $v$ denote a nonzero vector in $E^*_2V$ which is orthogonal to 
$s_2$.
Then the nonvanishing vectors among \begin{equation}\label{EIMVBASIS1}
E_1v,E_2v,\ldots,E_{D-1}v
\end{equation}
form an orthogonal basis for Mv.
\end{lemma}
{\it Proof.} Recall $E_0,E_1,\ldots,E_D$ form a basis for $M$.
We assume $v$ is orthogonal to $s_2$ so $E_0v=0$, 
$E_{D}v=0$.
Now apparently the vectors in (\ref{EIMVBASIS1}) span $Mv$. The 
vectors in (\ref{EIMVBASIS1}) are mutually orthogonal by 
(\ref{EIVDECOM}) and the result follows. \hfill $\Box$\\

\noindent
Referring to Lemma
\ref{EIMVBASIS}, in this section we determine which of $E_1v, E_2v, 
\ldots, E_{D-1}v$ are zero. To do this, we compute
the square norms of these vectors.

\begin{theorem}
\label{EIV}
With reference to Definition \ref{A},
let $v$ denote a nonzero vector in $E^*_2V$ which is orthogonal to 
$s_2$.
Assume $v$ is an eigenvector for
$E^*_2A_{2}E^*_2$; let $\eta $ denote the corresponding
eigenvalue.
\begin{description}
\item[(i)]
Assume $\eta \not=-1$. Then \begin{equation}\label{EIV1}
\|E_iv\|^2
= \frac{m_i(\theta_i-k)(\theta_{i}+k)(\theta_i^{2}-\psi)}
{|X|k b_{1}(\psi - b_{2})}
\|v\|^2
\qquad \qquad (0 \le i \le D),
\end{equation} where \begin{equation}  \label{psidef} \psi = b_{2}\left(1- 
\frac{b_{3}}{1+\eta}\right). \end{equation}
We remark the denominator in (\ref{EIV1}) is nonzero by the form of 
(\ref{psidef}).  
\item[(ii)]
Assume $\eta =-1$. Then \begin{equation}\label{EIV2}
\|E_iv\|^2
= \frac{m_i(k - \theta_i)(k + \theta_i)}{|X|kb_1} \|v\|^2
\qquad \qquad (0 \le i \le D).
\end{equation}
\end{description}
\end{theorem}
{\it Proof.} First assume $i \in \{0,D\}$.  Then $\theta_{i} \in \{-k, 
k\}$.  We assume $v$ is orthogonal to $s_2$ so $E_0v=0$, 
$E_{D}v=0$.  Thus the equations in (\ref{EIV1}), (\ref{EIV2}) hold 
since both sides of the equations are zero.  Now assume $1 \le i \le 
D-1$.  
Observe $v = E_2^*v$ by construction so $E_iv=E_iE^*_2v$.
We may now argue
\begin{eqnarray} 
\|E_iv\|^2
&=& (E_iE_2^*v)^t \overline{E_iE_2^*v} \qquad \qquad \qquad \;
{\rm by}\;(\ref{DOTDEF})\nonumber\\
{} &=& 
v^tE_2^{*t}E_i^t\overline{E}_i\overline{E}_2^*\overline{v}\nonumber\\
{} &=& v^tE_2^*E_i^2E_2^*\overline{v} \qquad \qquad \qquad \qquad
{\rm by}\;(eiii), (eiv), (dii), 
(diii)\nonumber\\
{} &=& v^tE_2^*E_iE_2^*\overline{v}
\qquad \qquad \qquad \quad \quad {\rm by}\;(ev).
\label{EIV3}
\end{eqnarray}
To evaluate (\ref{EIV3})
we apply Lemma
\ref{DEEDE}. To do this we make some comments.
Let
$\theta^*_0, \theta^*_1, \ldots, \theta^*_D$ denote the dual
eigenvalue sequence for $E_i$. We assume $v$ is orthogonal 
to $s_2$
so $Jv=0$. We already mentioned $E^*_2v=v$.
By assumption $E^*_2A_{2}¥E^*_2v=\eta v$.
Since each of $J$, $E^*_2$, $A_{2}¥$ has real entries and since $\eta $ is
real, we see
$J\overline{v} = 0$,
$E^*_2\overline{v} = {\overline{v}}$, and
$E^*_2A_{2}¥E^*_2\overline{v} = \eta \overline{v}$.
Evaluating (\ref{EIV3}) using Lemma \ref{DEEDE}
and our above comments, we find $|X|\,\|E_iv\|^2$ is equal to 
$\|v\|^2 $ times
\begin{equation}
\theta_0^*- \theta_4^*
+ \eta (\theta_2^* - \theta_4^*).
\label{EIV4}
\end{equation}
Observe (\ref{EIV4}) is equal to $\theta_{0}^{*}-\theta_{2}^{*}$ 
times 
\begin{equation} \label{EIV5}
  1 + (1+\eta)\frac{\theta_2^* - \theta_4^*}{\theta_{0}^{*}-\theta_{2}^{*}}.
\end{equation}
  Using (\ref{recth}) and recalling $\theta_{0}^{*} =m_{i}$, we find 
  $\theta_{0}^{*}-\theta_{2}^{*} = 
  m_{i}(k-\theta_{i})(k+\theta_{i})/(kb_{1})$.  First assume $\eta 
  =-1$.  Then (\ref{EIV5}) is equal to $1$ and (\ref{EIV2}) follows.  
  Next assume $\eta \not= -1$.  Evaluating (\ref{EIV5}) using 
  Corollary \ref{thfrac} and (\ref{psidef}), we routinely verify 
  that (\ref{EIV1}) holds.
\hfill $\Box$\\

\bigskip
\noindent
With reference to Definition \ref{SUBGRAPH}, 
our next goal is to get upper and lower bounds for the local 
eigenvalues $\eta_{i}$ $(k+1 \le i \le k_{2})$.  We will use the 
following notation.

\begin{definition}\label{TYPE}  \rm
With reference to Definition \ref{A},
for all $z \in \C \cup \infty$ 
we define
$${{\tilde z} = \cases{
-1-\frac{b_2 b_{3}}{z^{2}-b_{2}¥}, & if $z\not=\infty, \;z^{2} \ne {b_{2}}$ \cr
\infty, & if $z^{2} = {b_{2}}$ \cr
-1, &if $z=\infty $. \cr
}}$$
By Lemma
\ref{bipeigs}
neither of $\theta_1^{2}, \theta_d^{2}$ is equal to ${b_{2}}$,
so
\begin{equation}  \label{tildeth1d}
{\tilde \theta}_1=-1-b_2 b_{3}(\theta_1^{2}-b_{2})^{-1}, \qquad \qquad 
{\tilde \theta}_d=-1-b_2 b_{3}(\theta_d^{2}-b_{2})^{-1}. 
\end{equation}
By the data in
Lemma \ref{bipeigs} we have
${\tilde \theta}_1 < -1$. Moreover
${\tilde \theta}_d > b_3-1$ if $D$ is odd and  
${\tilde \theta}_d = b_3-1$ if $D$ is even.
In either case 
${\tilde \theta}_d \geq 0$.
\end{definition}

\begin{theorem}\label{TYPEINEQ}
With reference to Definitions \ref{A}
and \ref{SUBGRAPH},
we have
${\tilde \theta}_1 \leq \eta_i \leq {\tilde \theta}_d$ for
$k+1 \le i \le k_{2}$.
\end{theorem}
{\it Proof.} Let the integer $i$ be given, and abbreviate 
$\eta=\eta_i$.
Let $v$ denote a nonzero vector in $U_\eta$.
First suppose $\eta < {\tilde \theta}_1$. By Definition
\ref{TYPE} and Lemma \ref{bipeigs} we find
\begin{equation}  \label{tildecond1}
\eta +1 < \frac{b_{2}b_{3}}{b_{2}-\theta_{1}^{2}} < 0.
\end{equation}
Using (\ref{psidef}) and (\ref{tildecond1}), we find $\psi > b_{2}$
and $\psi < 
\theta_{1}^{2}$.  Now ${\|E_1v\|^2}<0$ by Theorem
\ref{EIV}, a contradiction. The inequality $\eta_{i} \le {\tilde 
\theta}_{d}$ is proven similarly.
\hfill $\Box$\\

\noindent Referring to Lemma
\ref{EIMVBASIS}, we now determine which of $E_1v, E_2v, \ldots, E_{D-1}v$ 
are zero.

\begin{lemma} \label{EJV}
With reference to Definition \ref{A},
let $v$ denote a nonzero vector in $U$. Then (i)--(vi) hold below. 
\begin{description}
\item[(i)] $E_0v=0$ and $E_{D}v=0$.
\item[(ii)] For $1 \le i \le D-1$, $\; E_{i}v \not= 0$ 
provided $i$ is not among $1, d, D-d, D-1$. 
\item[(iii)] $E_1v=0$  if and only if $v \in 
U_{{\tilde \theta}_1}$.
\item[(iv)] $E_{D-1}v=0$ if and only if $v \in 
U_{{\tilde \theta}_1}$.
\item[(v)] $E_dv=0$ if and only if $v \in
U_{{\tilde \theta}_d}$.
\item[(vi)] $E_{D-d}v=0$ if and only if $v \in
U_{{\tilde \theta}_d}$.
\end{description}
\end{lemma}
{\it Proof.}  By Definition \ref{U=V0+Y}, $U$ is orthogonal to 
$V_{0}$; in particular $v$ is orthogonal to $s_{2}$.  
Thus $E_0v$, $E_{D}v$ are zero.
Now suppose there exists an integer $n$ $(1 \leq n 
\leq D-1)$
such that $E_nv=0$.
We show $n$ is among $1, d, D-d, D-1$, and that
$v\in U_{{\tilde \theta}_n}$.
We claim
$v$ is an eigenvector for $E^*_2A_{2}E^*_2$. To see this,
in (\ref{DEEDE1}) set $E=E_n$ and
apply both sides to $v$. Using $v=E_2^*v$ and $Jv=0$ we find 
$$0 = 
(\theta_0^*-\theta_4^*)v + (\theta_2^*-\theta_4^*)E^*_2A_{2}E^*_2v,$$
where $\theta^*_0, \theta^*_2, \theta^*_4$ are dual eigenvalues for 
$E_n$.  Observe $\theta_2^* \ne \theta_4^*$; otherwise $\theta_{0}^{*} = 
\theta_{4}^{*}$ by the above line, forcing $\theta_{0}^{*} = 
\theta_{2}^{*}$ and contradicting Lemma \ref{PITHETA*}.
Apparently $v$ is an eigenvector for $E^*_2A_{2}E^*_2$, as claimed. Let 
$\eta $ denote the corresponding eigenvalue.
By assumption $E_nv=0$ so $\Vert E_nv\Vert^2 =0$.
Applying Theorem \ref{EIV} we find $\theta_n^2 = \psi$ and
thus $ \eta={\tilde \theta}_n$. 
Thus 
$n$ is among $1, d, D-d, D-1$ by Theorem \ref{TYPEINEQ}, 
Lemma \ref{bipequiv}(v), and the fact that $\theta_{0} > 
\theta_{1} > \cdots > \theta_{D}$.  
 Apparently
$v \in U_{{\tilde \theta}_n}$.
To finish the proof, suppose $n=1$ or $n=d$ and assume $v \in
U_{{\tilde \theta}_n}$.
We show $E_nv=0$, $E_{D-n}v=0$.
Observe $v$ is an eigenvector for $E^*_2A_{2}E^*_2$ with eigenvalue 
${\tilde \theta}_n$.
Applying Theorem
\ref{EIV} we find $\Vert E_nv\Vert^2 =0$
so $E_nv=0$, as desired.  Observe $E_{D-n}v=0$ since 
$\Vert E_{D-n}v\Vert^2 = \Vert E_nv\Vert^2$ by (\ref{EIV1}).  
The result follows.
\hfill $\Box$\\

\begin{corollary} \label{dimj}
With reference to Definition \ref{A},
let $v$ denote a nonzero vector in $U$.
Then (i)--(iv) hold below.
\begin{description}
\item[(i)] If $v \in
U_{{\tilde \theta}_1}$ then
$Mv$ has dimension $D-3$.
\item[(ii)] If $v \in
U_{{\tilde \theta}_d}$ and $D$ is odd, then
$Mv$ has dimension $D-3$.
\item[(iii)] If $v \in
U_{{\tilde \theta}_d}$ and $D$ is even, then
$Mv$ has dimension $D-2$.
\item[(iv)]
If $v \notin
U_{{\tilde \theta}_1}$
and
$v \notin
U_{{\tilde \theta}_d}$ then
$Mv$ has dimension $D-1$.
\end{description}
\end{corollary}
{\it Proof.} Combine
Lemmas \ref{EIMVBASIS}
and \ref{EJV} and observe the integers $d$, $D-d$ are distinct 
precisely when $D$ is odd.
\hfill $\Box$\\

\noindent
The following equations will be useful.  

\begin{lemma} \label{multeqs}
With reference to Definitions \ref{A} and \ref{phidef}, the following 
(i)--(iii) hold.
\begin{description}
\item[(i)] $\displaystyle{k + \sum_{\eta \in \Phi} 
\mbox{mult}_{\eta} =k_{2}}$.
\item[(ii)] $\displaystyle{p_{22}^{2} +(k-1)(b_{3}-1) +\sum_{\eta \in \Phi} 
\eta \, \mbox{mult}_{\eta}=0}$.
\item[(iii)]  $\displaystyle{(p_{22}^{2})^{2} +(k-1)(b_{3}-1)^{2} +
\sum_{\eta \in \Phi} 
\eta^{2} \mbox{mult}_{\eta}=k_{2}p_{22}^{2}}$.
\end{description}
\end{lemma}
{\it Proof.}  (i) There are $k_{2}-k$ elements in the sequence 
$\eta_{k+1}, \eta_{k+2}, \ldots, \eta_{k_{2}}$. \\
(ii)  Recall the matrix $\breve{A}$ from Definition \ref{SUBGRAPH}.
Each diagonal entry of $\breve{A}$ is zero, so the trace of 
$\breve{A}$ is zero.  Recall $\eta_{1}, \eta_{2}, \ldots, \eta_{k_{2}}$ 
are the eigenvalues of $\breve{A}$, so $\sum_{i=1}^{k_{2}}\eta_{i} 
=0$.  By this and since $\eta_{1} =p_{22}^{2}$, $\eta_{i} =b_{3}-1 \; 
(2 \le i \le k)$, we have the desired result.  \\
(iii)  Recall $\Gamma_{2}^{2}$ is regular with valency $p_{22}^{2}$, 
so each diagonal entry of $\breve{A}^{2}$ is $p_{22}^{2}$.  
Apparently the trace of $\breve{A}^{2}$ is $k_{2}p_{22}^{2}$, so 
$\sum_{i=1}^{k_{2}} \eta_{i}^{2} =k_{2}p_{22}^{2}$.  By this and 
since $\eta_{1}=p_{22}^{2}, \; \eta_{i}=b_{3}-1 \; (2 \le i \le k)$, 
we have the desired result.
\hfill $\Box$\\

\begin{definition}
\label{TM}  \rm
With reference to Definition \ref{A}, let
$W$ denote a thin irreducible $T$-module with endpoint 2.
Observe $E^*_2W$ is a $1$-dimensional eigenspace for
$E^*_2A_{2}¥E^*_2$; let $\eta $ denote the corresponding eigenvalue.
We observe $E^*_2W$ is contained in $E^*_2V$ and is orthogonal
to any irreducible $T$-module with endpoint 0 or 1, so
$E^*_2W\subseteq U_\eta$.
Apparently $U_\eta \not=0$ so $\eta $
is one of $\eta_{k+1}, \eta_{k+2},\ldots, \eta_{k_{2}}$.
We have ${\tilde \theta}_1 \leq
\eta \leq {\tilde \theta}_d$ by Theorem \ref{TYPEINEQ}. We refer to 
$\eta $ as the {\em local eigenvalue} of $W$.
\end{definition}

\noindent With reference to Definition \ref{A}, let $W$ denote a thin 
irreducible $T$-module with endpoint 2 and local eigenvalue
$\eta $. In order to describe $W$ we distinguish four cases: (i) $D$ 
is odd, and $\eta 
={\tilde \theta}_1 $ or $\eta ={\tilde \theta}_d $; (ii) $D$ is even 
and $\eta = {\tilde \theta}_{1}$; (iii) $D$ is even 
and $\eta = {\tilde \theta}_{d}$; (iv) $ {\tilde \theta}_1 < \eta < {\tilde 
\theta}_d$.  
We investigate cases (i), (ii) in the present paper.
We will investigate the remaining cases in a future paper.


\section{The spaces $U_{{\tilde \theta}_1}$ and $U_{{\tilde 
\theta}_d}$ }

\noindent We state our main goal for this section.
With reference to Definition \ref{A}, choose $n \in \{1,d\}$ if $D$
is odd, and let $n=1$ if $D$ is even.  Define $\eta = {\tilde \theta}_n$.
We show that for all nonzero $v\in U_\eta $
the space $Mv$ is a thin
irreducible $T$-module with endpoint 2 and local eigenvalue $\eta$.

\begin{lemma}
\label{C1}
With reference to Definition \ref{A}, let
$v$ denote a vector in $E^*_2V$.
Then
$$
E^*_iA_jv = 0 \qquad \hbox{if} \qquad |i-j|>2 \qquad \qquad (0 \leq 
i,j\leq D), 
$$ and
$$E^*_iA_jv = 0 \qquad \hbox{if} \qquad i+j \hbox{ is odd} \qquad \qquad (0 \leq 
i,j\leq D). $$
\end{lemma}
{\it Proof.}
Let $i,j$ be given and observe $E^*_2v=v$ so $E^*_iA_jv = E^*_iA_jE^{*}_{2}v$.
The result now follows from (\ref{REL2}).
\hfill $\Box$\\

\begin{lemma}
\label{C2}
With reference to Definition \ref{A}, let
$v$ denote a vector in $E^*_2V$ which is
orthogonal to $s_2$.
Then 
$$
\sum_{j=0}^D E^*_iA_jv = 0 \qquad \qquad (0 \leq i \leq D).
$$
\end{lemma}
{\it Proof.} Observe $Jv=0$ so $E^*_iJv=0$. Eliminate $J$ in this
expression using (aii) to get
the result.
\hfill $\Box$\\

\begin{lemma}
\label{Cspl}
With reference to Definition \ref{A}, choose $n \in \{1,d\}$ if $D$
is odd, and let $n=1$ if $D$ is even.
Define $\eta = {\tilde \theta}_n$. Then for all
$v\in U_\eta$ we have 
$$
\sum_{j=0}^D \theta^*_jE^*_iA_jv = 0 \qquad \qquad (0 \leq i \leq D),
$$
where $\theta^*_0, \theta^*_1, \ldots, \theta^*_D$
denotes the dual eigenvalue sequence for $\theta_n$.
\end{lemma}
{\it Proof.} Observe $E_nv=0$
by Lemma \ref{EJV} so $E^*_iE_nv=0$. Eliminate $E_n$ in this
expression using
(\ref{PRIMIDEM})
to get
the result.
\hfill $\Box$\\

\begin{lemma}\label{IDEN1}
With reference to Definition \ref{A}, choose $n \in \{1,d\}$ if $D$
is odd, and let $n=1$ if $D$ is even.
Define $\eta = {\tilde \theta}_n$. Then for all $v\in U_\eta$ we 
have 
\begin{eqnarray} 
E_i^*A_iv &=& \frac{\theta_{i-2}^*-\theta_{i+2}^*}{\theta_{i+2}^{*}-\theta_{i}^{*}}
 E_i^*A_{i-2}v \qquad \qquad 
(2 \le i \le D-2),
\label{IDEN1A}\\
E_i^*A_{i+2}v &=& \frac{\theta_{i-2}^*-\theta_{i}^*}{\theta_{i}^{*}-\theta_{i+2}^{*}}
 E_i^*A_{i-2}v \qquad \qquad
(2 \le i \le D-2), \label{IDEN1B}
\end{eqnarray}
where $\theta^*_0, \theta^*_1, \ldots, \theta^*_D$
denotes the dual eigenvalue sequence for $\theta_n$.  
Moreover \begin{equation}\label{IDEN1C}
E_0^*A_iv = 0, \quad E_1^*A_iv = 0, \quad E_{D-1}^*A_iv = 0, \quad E_D^*A_iv = 0
\qquad \qquad (0 \leq i \leq D).
\end{equation}
We note the denominators in (\ref{IDEN1A}), (\ref{IDEN1B}) are 
 nonzero by Lemmas \ref{PITHETA*}, \ref{PITHETA}, and \ref{PITHETAD}.
\end{lemma}
{\it Proof.} By Lemmas \ref{PITHETA*}, \ref{PITHETA}, and \ref{PITHETAD}, we
find $\theta^*_{i-2} \not = \theta^*_i$ for $2 \le i \le D$.
Solving the equations in Lemma \ref{C2} and
Lemma \ref{Cspl} using Lemma \ref{C1},
we routinely obtain (\ref{IDEN1A})--(\ref{IDEN1C}).
\hfill $\Box$\\

\begin{lemma} \label{C3}
With reference to Definition \ref{A}, let
$v$ denote a vector in $E^*_2V$ which is
orthogonal to $s_2$. Let the polynomials
$p_0, p_1, \ldots, p_D$ be from (\ref{PIPOLY}).
Then \begin{equation}\label{C3e}
p_i(A)v = E^*_{i+2}A_iv-E^*_iA_{i+2}v \qquad \qquad (0 \leq i \leq 
D-2). \end{equation}
Moreover $p_{D-1}(A)v=0, \; p_D(A)v=0$. \end{lemma}
{\it Proof.} For $0 \leq i \leq D-2$ we have \begin{eqnarray}
p_i(A)v &=& \sum_{{0 \le s \le i}\atop{ i-s {\tiny \mbox{ 
even}}}}A_{s}v
\nonumber
\\
&=& (E^*_0+ E^*_1+\cdots + E^*_D)\sum_{{0 \le s \le i}\atop{ i-s {\tiny \mbox{ 
even}}}}A_{s}v
\nonumber
\\
&=& \sum E^*_rA_sv, \label{C3e2}
\end{eqnarray}
where the final sum is over all integers $r,s$ such that $0 \leq r \leq 
D$, $0 \leq s \leq i$, and $i-s$ is even. Cancelling terms in 
(\ref{C3e2}) using
Lemmas \ref{C1} and \ref{C2} we obtain (\ref{C3e}).  Using 
(\ref{PDAeqs}), we may solve for each of $p_{D-1}(A), \; p_{D}(A)$ as a 
linear combination of $J$ and $J'$.  Recall 
$Jv=0$, $J'v=0$, and thus $p_{D-1}(A)v=0, \; p_D(A)v=0$.
\hfill $\Box$\\

\begin{theorem} \label{T1}
With reference to Definition \ref{A}, choose $n \in \{1,d\}$ if $D$
is odd, and let $n=1$ if $D$ is even.
Define $\eta = {\tilde \theta}_n$. Then for all $v\in U_\eta$ we 
have \begin{equation}\label{T1a}
\label{IDEN4A}
E_{i+2}^{*}A_{i}v = 
\sum_{{0 \le h \le i}\atop{ i-h {\tiny \mbox{ even}}}}
\frac{\theta^*_h-\theta^*_{h+2}}{\theta_{i}^{*}-\theta_{i+2}^{*}}
p_h(A)v
\qquad \qquad (0 \leq i \leq D-2),
\end{equation}
where $\theta^*_0, \theta^*_1, \ldots, \theta^*_D$
denotes the dual eigenvalue sequence for $\theta_n$.
Moreover each side of (\ref{IDEN4A}) is zero for $i=D-3$, $i=D-2$.  
We note the denominators in (\ref{IDEN4A}) are 
 nonzero by Lemmas \ref{PITHETA*}, \ref{PITHETA}, and \ref{PITHETAD}.
\end{theorem}
{\it Proof.}
To verify
(\ref{IDEN4A}),
in the expression on the right
eliminate $p_h(A)v$
using (\ref{C3e}), and simplify the result using
(\ref{IDEN1B}),
(\ref{IDEN1C}). We now have (\ref{IDEN4A}).  For $i=D-3$, $i=D-2$, 
each side of (\ref{IDEN4A}) is zero by (\ref{IDEN1C}).  
\hfill $\Box$\\

\begin{theorem}  \label{EiGi}
With reference to Definition \ref{A}, choose $n \in \{1,d\}$ if $D$
is odd, and let $n=1$ if $D$ is even.
Define 
$\eta = {\tilde \theta}_n$.  Let $\theta=\theta_{n}$, and let 
$g_{0}, g_{1}, \ldots, g_{D-2}$ denote the associated polynomials from 
Definition \ref{GIPOLY}.  
Then for all $v \in U_\eta $ and for $0 \le i \le D-2$,
\begin{equation}  \label{gia}
  E_{i+2}^{*}A_{i}v = g_{i}(A)v.
\end{equation}
Moreover, each side of (\ref{gia}) is zero for $i=D-3$, $i=D-2$.  
\end{theorem}
{\it Proof.}  Choose an integer $h \; (0 \le h \le D-2)$.
Using Lemma \ref{PITHETA*} and (\ref{KI}), we find 
$$  \frac{k_{i}b_{i}b_{i+1}}{k_{h}b_{h}b_{h+1}} 
\frac{p_{h}(\theta)}{p_{i}(\theta)} = 
\frac{\theta_{h}^{*}-\theta_{h+2}^{*}}{\theta_{i}^{*}-\theta_{i+2}^{*}}, $$
where denominators are nonzero by Lemmas
\ref{PITHETA}, \ref{PITHETAD}.  The result now follows by Definition 
\ref{GIPOLY} and Theorem \ref{T1}.  For $i=D-3$, $i=D-2$ both sides 
of (\ref{gia}) are zero by Theorem \ref{T1}.  
\hfill $\Box$\\


\begin{lemma} \label{mvbsis}
With reference to Definition \ref{A}, choose $n \in \{1,d\}$ if $D$
is odd, and let $n=1$ if $D$ is even.
Define $\eta = {\tilde \theta}_n$. 
Then for all nonzero $v \in U_\eta, $ 
the vectors
$
E^*_{i+2}A_iv \; (0 \leq i \leq D-4)
$
form a basis for $Mv$.  
\end{lemma}
{\it Proof.}
By Corollary \ref{dimj} the dimension of
$Mv$ is $D-3$. By this and since $A$ generates $M$ we find $v, Av, 
A^2v, \ldots, A^{D-4}v$
form a basis for $Mv$. For $0 \leq i \leq D-2$
let the polynomial $g_i$ be as in Definition \ref{GIPOLY} (with 
$\theta = \theta_{n}$).  Recall 
$g_i$ has degree $i$.
Apparently the vectors $g_i(A)v$ $(0 \leq i \leq D-4)$
form a basis for $Mv$.   By Theorem \ref{EiGi} we have $g_{i}(A)v =
E_{i+2}^{*}A_{i}v$ for $0 \le i \le D-4$.  The result follows. 
\hfill $\Box$\\

\begin{theorem}
\label{Gonew4}
With reference to Definition \ref{A}, choose $n \in \{1,d\}$ if $D$
is odd, and let $n=1$ if $D$ is even.  Define $\eta = {\tilde \theta}_n$.
Then for all nonzero $v \in U_\eta $
the space $Mv$ is a thin
irreducible $T$-module with endpoint 2 and local eigenvalue $\eta$.
\end{theorem}
{\it Proof.}
We first show $Mv$ is a $T$-module. It is clear
$Mv$ is closed under $M$.
By Lemma
\ref{mvbsis}
and (div) we find
$Mv$ is closed under $M^*$. Recall $M$ and $M^*$ generate $T$ so $Mv$ 
is a $T$-module.
We show $Mv$ is irreducible. From Lemma
\ref{mvbsis}
we find $v$ is a basis for $E^*_2Mv$.
In particular
$E^*_2Mv$ has dimension 1. Since $Mv$ is a $T$-module it is a
direct sum of irreducible $T$-modules. It follows there exists
an irreducible $T$-module $W'$ such that $W'\subseteq Mv$ and
such that $E^*_2W'\not=0$. We show $W'=Mv$. Observe
$E^*_2W' \subseteq E^*_2Mv$, and we mentioned
$E^*_2Mv$ has dimension 1, so
$E^*_2W' =E^*_2Mv$. Now apparently $v \in E^*_2W'$. Observe $W'$ is
$M$-invariant, so $Mv\subseteq W'$, and it follows $W'=Mv$.
In particular $Mv$ is irreducible. 
From Lemma \ref{mvbsis}
we find $E^*_iMv$ is 0
for $i\in \lbrace 0,1,D-1,D \rbrace $ and has dimension 1 for $2 \leq i 
\leq D-2$. Apparently
$Mv$ is thin with endpoint 2.
We mentioned $v$ is a basis for $E^*_2Mv$. From the construction
$v\in U_\eta
$ so
$Mv$ has local eigenvalue $\eta$.  
\hfill $\Box$\\


\section{The thin irreducible $T$-modules with endpoint 2
and local eigenvalue ${\tilde \theta}_1$ or ${\tilde \theta}_d$ }

\noindent With reference to Definition \ref{A}, choose $n \in \{1,d\}$ if $D$
is odd, and let $n=1$ if $D$ is even.
We now
describe the thin irreducible $T$-modules
with endpoint 2
and local eigenvalue
${\tilde \theta}_n$.

\begin{theorem} \label{T5} With reference to Definition \ref{A}, 
choose $n \in \{1,d\}$ if $D$
is odd, and let $n=1$ if $D$ is even.  Let
$W$ denote a thin irreducible $T$-module with endpoint 2
and local eigenvalue
${\tilde \theta}_n$.
Let $v$ denote
a nonzero vector in $E^*_2W$. Then $W=Mv$. The vectors 
\begin{equation}\label{T5a}
E_iv \qquad (1 \leq i \leq D-1, \;\;i\not=n, \;\; i\not= D-n)
\end{equation}
form a basis for $W$ and $E_0v=0$, $E_nv=0$, $E_{D-n}v=0$, $E_{D}v=0$.
\end{theorem}
{\it Proof.}
Observe $W$ is $M$-invariant and $v\in W$ so
$Mv\subseteq W$. Observe $v \in U_{{\tilde \theta}_n}$ by Definition
\ref{TM}; combining this with Theorem
\ref{Gonew4} we find $Mv$ is a $T$-module. Now
$W=Mv$ by the irreducibility of $W$. We mentioned
$v \in U_{{\tilde \theta}_n}$; by this and
Lemma \ref{EJV} we find each of the vectors in (\ref{T5a}) are 
nonzero. Moreover $E_0v=0$, $E_nv=0$, $E_{D-n}v=0$, and $E_{D}v=0$. 
Applying Lemma
\ref{EIMVBASIS} we find the vectors in (\ref{T5a}) form a basis for 
$Mv$.
\hfill $\Box$\\

\begin{theorem} \label{T6}
With reference to Definition \ref{A}, choose $n \in \{1,d\}$ if $D$
is odd, and let $n=1$ if $D$ is even.  Let
$W$ denote a thin irreducible $T$-module with endpoint 2
and local eigenvalue
${\tilde \theta}_n$.
The vectors in (\ref{T5a}) are mutually orthogonal and 
\begin{equation}\label{T6a}
\|E_iv\|^2 = \frac{m_i(\theta_i^{2}-k^{2})(\theta_i^{2}-\theta_{n}^{2})}
{|X|k b_{1}(\theta_{n}^{2} - b_{2})}
\|v\|^2
\qquad \qquad (1 \le i \le D-1, \;\;i\not=n, \;\; i\not= D-n),
\end{equation}
where the scalar $m_i$ denotes the multiplicity of $\theta_i$.  We 
remark the denominator in (\ref{T6a}) is nonzero by Lemma 
\ref{bipeigs}.  
\end{theorem}
{\it Proof.} The vectors in (\ref{T5a}) are mutually orthogonal by
(\ref{EIVDECOM}). To obtain (\ref{T6a}) we apply
Theorem
\ref{EIV}. Set $\eta = {\tilde \theta}_n$ and
observe $\eta \not=-1$ by
Definition
\ref{TYPE}.
Now
(\ref{EIV1}) holds and (\ref{T6a}) follows.
\hfill $\Box$\\

\begin{theorem} \label{T2}
With reference to Definition \ref{A}, choose $n \in \{1,d\}$ if $D$
is odd, and let $n=1$ if $D$ is even.  Let
$W$ denote a thin irreducible $T$-module with endpoint 2
and local eigenvalue ${\tilde \theta}_n$.
Let $v$ denote
a nonzero vector in $E^*_2W$. 
Then 
the vectors
$
E^*_{i+2}A_iv \; (0 \leq i \leq D-4)
$
form a basis for $W$. 
\end{theorem}
{\it Proof.}  Observe $v \in U_{{\tilde \theta}_n}$ by
Definition
\ref{TM} and $W=Mv$ by Theorem
\ref{T5}. The result now follows in view
of Lemma
\ref{mvbsis}.
\hfill $\Box$\\

\begin{theorem} \label{T7}
With reference to Definition \ref{A}, choose $n \in \{1,d\}$ if $D$
is odd, and let $n=1$ if $D$ is even.  Let
$W$ denote a thin irreducible $T$-module with endpoint 2
and local eigenvalue
${\tilde \theta}_n$.  Abbreviate $\theta=\theta_{n}$, and 
let $g_{0}, g_{1}, \ldots, g_{D-2}$ denote the associated polynomials 
from Definition \ref{GIPOLY}.  
Let $v$ denote a nonzero vector in $E^*_2W$.
Then for $0 \leq i \leq D-4$ we have \begin{equation}\label{T7a}
E^*_{i+2}A_iv = \sum_{{1 \le j \le D-1}\atop{j \not=n, \, j \not=D-n}}
g_i(\theta_j) 
E_jv.
\end{equation}
\end{theorem}
{\it Proof.} By Theorem
\ref {EiGi} we have
$E^*_{i+2}A_iv=g_i(A)v$. In this equation, multiply $g_i(A)v$
on the left by $I$, expand using (eii), and simplify the result 
using $AE_j=\theta_jE_j$
$(0 \leq j \leq D)$.  Observe $E_{0}v=0$, $E_{n}v=0$, $E_{D-n}v=0$, 
$E_{D}v=0$ by Theorem \ref{T5}.  \hfill $\Box$\\

\begin{theorem} \label{T3}
With reference to Definition \ref{A}, choose $n \in \{1,d\}$ if $D$
is odd, and let $n=1$ if $D$ is even.  Let
$W$ denote a thin irreducible $T$-module with endpoint 2
and local eigenvalue ${\tilde \theta}_n$.
Let $v$ denote a nonzero vector in $E_2^*W$.
The vectors $E_{i+2}^{*}A_{i}v \; (0 \le i \le D-4)$
 are mutually orthogonal and 
\begin{equation}\label{T3a}
\|E^*_{i+2}A_iv\|^2 = 
\frac{k_{i}b_{i}b_{i+1}c_{i+1}c_{i+2}}{kb_{1}(\theta_{n}^{2} - b_{2})}
\frac{p_{i+2}(\theta_{n})}{p_{i}(\theta_{n})} \| v \|^{2}  \qquad (0 
\le i \le D-4),
\end{equation}
where the polynomials $p_{i}$ are as in (\ref{PIPOLY}).  We 
remark the denominators in (\ref{T3a}) are nonzero by Lemmas 
\ref{bipeigs}, \ref{PITHETA}, and \ref{PITHETAD}.  
\end{theorem}
{\it Proof.}  The vectors $E_{i+2}^{*}A_{i}v \; (0 \le i \le D-4)$
are mutually orthogonal by
(\ref{EISORTHO}). To verify (\ref{T3a}), in the left-hand side 
eliminate $E^*_{i+2}A_iv$
using (\ref{T7a}), and evaluate the result using
Theorem
\ref{GIGJ} and Theorem \ref{T6}.
\hfill $\Box$\\

\begin{theorem} \label{T4}
With reference to Definition \ref{A}, 
choose $n \in \{1,d\}$ if $D$
is odd, and let $n=1$ if $D$ is even.  Let
$W$ denote a thin irreducible $T$-module with endpoint 2
and local eigenvalue ${\tilde \theta}_n$. 
With respect to the basis for $W$ given in Theorem \ref{T2} the matrix 
representing $A$ is
$$
\left(\begin{array}{cccccc}
0 & \omega_1 & & & & {\bf 0}\\
c_1 & 0 & \omega_2 & & & \\
& c_2 & \cdot & \cdot & & \\
& & \cdot & \cdot & \cdot& \\
& & & \cdot & \cdot & \omega_{D-4} \\
{\bf 0} & & & & c_{D-4} & 0
\end{array} \right),
$$
where 
\begin{equation}  \label{omega2}
  \omega_{i} = \frac{b_{i+1}c_{i+2}}{c_{i}} \frac{p_{i-1}(\theta_{n}) 
  p_{i+2}(\theta_{n})}{p_{i}(\theta_{n}) p_{i+1}(\theta_{n})} \qquad 
  \qquad (1 \le i \le D-4).
\end{equation}
We remark the denominator in (\ref{omega2}) is nonzero by Lemmas 
\ref{PITHETA} and \ref{PITHETAD}.  
\end{theorem}
{\it Proof.}  Let $\theta=\theta_{n}$, and 
let $g_{0}, g_{1}, \ldots, g_{D-2}$ denote the associated polynomials 
from Definition \ref{GIPOLY}.  Setting $\lambda = A$ in (\ref{girecur})
and applying the result to $v$, we find
$$
Ag_i(A)v= c_{i+1}g_{i+1}(A)v +  \omega_i g_{i-1}(A)v
\qquad \quad (0 \leq i \leq D-4),
$$
where $g_{-1}=0$, $\omega_{0}=0$. The result follows in view of
Theorem
\ref{EiGi}.
\hfill $\Box$\\

\noindent In summary we have the following theorem.

\begin{theorem}
\label{SUM}
With reference to Definition \ref{A}, choose $n \in \{1,d\}$ if $D$
is odd, and let $n=1$ if $D$ is even.   Let
$W$ denote a thin irreducible $T$-module with endpoint 2
and local eigenvalue
${\tilde \theta}_n$.
Then $W$ has dimension $D-3$.
For $0 \leq i \leq D$, $E^*_iW$ is zero if $i\in \lbrace 0,1,D-1,D\rbrace$
and has dimension
1 if
$i\not\in \lbrace 0,1,D-1,D\rbrace$.
Moreover
$E_iW$ is zero if $i\in \lbrace 0,n, D-n, D\rbrace $ and has dimension
1 if $i\not\in \lbrace 0,n, D-n, D\rbrace $.
\end{theorem}
{\it Proof.}  The dimension of $W$ is equal to $D-3$ by Theorem
\ref{T2}.
Fix an integer $i$ $(0 \leq i \leq D)$.
From Theorem
\ref{T2} we find
$E^*_iW$ is zero if
$i\in \lbrace 0,1,D-1,D\rbrace$
and has dimension
1 if
$i\not\in \lbrace 0,1,D-1,D\rbrace$.
From Theorem
\ref{T5}
we find
$E_iW$ is zero if $i\in \lbrace 0,n,D-n,D\rbrace $ and has dimension
1 if $i\not\in \lbrace 0,n,D-n,D\rbrace $. \hfill $\Box$\\


\section{Some multiplicities  }  

With reference to Definition \ref{A}, choose $n \in \{1,d\}$ if $D$
is odd, and let $n=1$ if $D$ is even.   Let $W$ denote a thin
irreducible $T$-module with endpoint 2 and local eigenvalue
${\tilde \theta}_n$.
In this section we consider the multiplicity with which $W$ appears 
in the standard module $V$.

\begin{theorem} \label{T8}
With reference to Definition \ref{A}, choose $n \in \{1,d\}$ if $D$
is odd, and let $n=1$ if $D$ is even.   Let
$W$ denote a thin irreducible $T$-module with endpoint 2
and local eigenvalue ${\tilde \theta}_n$.
Let $W'$ denote an irreducible $T$-module. Then the following
(i), (ii) are equivalent.
\begin{description}
\item[(i)] $W $ and $W'$ are isomorphic as $T$-modules.
\item[(ii)]
$W'$ is thin with endpoint 2 and local eigenvalue ${\tilde \theta}_n$.
\end{description}
\end{theorem}
{\it Proof.} (i) $\Rightarrow $ (ii) Clear.

\noindent (ii) $\Rightarrow $ (i)
We display an isomorphism of $T$-modules from $W$ to $W^\prime $. 
Observe $E^*_2W$ and $E^*_2W'$ are both nonzero.
Let $v$ (resp. $v'$) denote a nonzero
vector in $E^*_2W$ (resp. in $E^*_2W'$). By Theorem
\ref{T2} the vectors
\begin{equation}
E^*_{i+2}A_iv \qquad \qquad (0 \leq i \leq D-4)
\label{b1}
\end{equation} form
a basis for $W$. Similarly the vectors
\begin{equation}
E^*_{i+2}A_iv' \qquad \qquad (0 \leq i \leq D-4)
\label{b1p}
\end{equation} form
a basis for $W'$.
Let $\sigma : W \rightarrow W'$ denote the isomorphism of
vector spaces that sends
$E^*_{i+2}A_iv$
to $E^*_{i+2}A_iv'$ for
$0 \leq i \leq D-4$. We show $\sigma $ is an isomorphism
of $T$-modules. By Theorem
\ref{T4}
the matrix representing $A$ with respect to the basis
(\ref{b1})
is equal to the matrix
representing $A$ with respect to the basis
(\ref{b1p}).
It follows $ \sigma A - A \sigma $ vanishes on $W$.
From the construction we find
that for $0 \leq h \leq D$,
the matrix representing $E^*_h$ with respect to the basis
(\ref{b1}) is equal to the matrix
representing $E^*_h$ with respect to the basis
(\ref{b1p}).
It follows $ \sigma E^*_h - E^*_h \sigma $ vanishes on $W$.
The algebra $T$ is generated by
$A,
E_0^*,E_1^*,\ldots,E_D^*$.
It follows $\sigma B - B \sigma $ vanishes
on $W$ for all $B \in T$. We now see $\sigma $ is an isomorphism
of $T$-modules from $W$ to $W^\prime $.  
\hfill $\Box$\\

\begin{lemma}
\label{T11}
With reference to Definition \ref{A}, choose $n \in \{1,d\}$ if $D$
is odd, and let $n=1$ if $D$ is even.  Define $\eta = {\tilde \theta}_n$.
Then
$$U_\eta = E^*_2H_\eta,$$
where $H_\eta$ denotes the subspace of $V$ spanned by all the thin 
irreducible $T$-modules with endpoint 2
and local eigenvalue $\eta$.
\end{lemma}
{\it Proof.}
We first show
$U_\eta \subseteq E^*_2H_\eta$.
Assume $U_\eta \not=0$; otherwise the result is trivial.
Let $v$ denote a nonzero vector in $U_\eta$.
By Theorem
\ref{Gonew4} we find $Mv$ is a thin irreducible
$T$-module with endpoint 2 and local eigenvalue $\eta$,
so $Mv \subseteq H_\eta$. Of course $v \in Mv$ so $v \in H_\eta$.
By the construction $v \in E^*_2V $ so $v=E^*_2v$. It follows
$v \in E^*_2H_\eta$. We have now shown
$U_\eta \subseteq E^*_2H_\eta$. Next we show
$U_\eta \supseteq E^*_2H_\eta$. To see this observe
$E^*_2H_\eta$ is spanned by the $E^*_2W$, where $W$ ranges over all 
thin irreducible $T$-modules with
endpoint 2 and local eigenvalue $\eta$. For all such
$W$ the space $E^*_2W$ is contained in $U_\eta$
by Definition
\ref{TM}.
It follows $U_\eta \supseteq E^*_2H_\eta$.
\hfill $\Box$\\

\begin{definition}
\label{MULTD}  \rm
With reference to Definition \ref{A}, and from our discussion
in Section 7,
the standard module $V$ can be
decomposed into an orthogonal direct sum of irreducible $T$-modules.
Let $W$ denote an irreducible $T$-module. By the {\em multiplicity
with which $W$ appears in $V$}, we mean the number of irreducible
$T$-modules in the above decomposition which are isomorphic to $W$. 
\end{definition}

\begin{definition}
\label{D11A}  \rm
With reference to Definition \ref{A}, choose $n \in \{1,d\}$ if $D$
is odd, and let $n=1$ if $D$ is even.  Define 
$\eta = {\tilde \theta}_n$.
We let $\mu_\eta$ denote the multiplicity with which
$W$ appears in $V$, where $W$ is a thin irreducible $T$-module
with endpoint 2 and local eigenvalue $\eta$.
If no such $W$ exists we set $\mu_\eta =0$.
\end{definition}

\begin{theorem} \label{T12}
With reference to Definition \ref{A}, choose $n \in \{1,d\}$ if $D$
is odd, and let $n=1$ if $D$ is even.  Define $\eta = {\tilde \theta}_n$.
Then the following scalars (i)--(iii) are equal:
\begin{description}
\item[(i)]
The scalar $\mu_\eta $ from Definition
\ref{D11A}.
\item[(ii)]
The dimension of $U_\eta $.
\item[(iii)]  The scalar $\mbox{mult}_{\eta}$ from Definition 
\ref{phidef}.
\end{description}
\end{theorem}
{\it Proof.} We mentioned below
Definition \ref{phidef} that the above scalars (ii), (iii) are 
equal. We now show the scalars (i), (ii) are equal.
In view of
Lemma \ref{T11} it suffices to show
the dimension of $E^*_2H_\eta $ is $\mu_\eta$. Observe $H_\eta$ is a 
$T$-module so it is an orthogonal direct sum
of irreducible $T$-modules. More precisely 
\begin{equation}\label{hometa}
H_\eta = W_1 + W_2 +\cdots + W_m \qquad \qquad {\rm (orthogonal\ 
direct\ sum)},
\end{equation}
where $m$ is a nonnegative integer, and
where $W_1,W_2,\ldots,W_m$ are thin irreducible $T$-modules with 
endpoint 2 and local eigenvalue $\eta$.
Apparently $m$ is equal to $\mu_\eta$.
We show $m$ is equal to the dimension of $E^*_2H_\eta $.
Applying $E_2^*$ to (\ref{hometa})
we find
\begin{equation}\label{hometa2}
E^*_2H_\eta = E_2^*W_1 + E_2^*W_2 +\cdots + E_2^*W_m \qquad \qquad 
{\rm (orthogonal\ direct\ sum)}.
\end{equation}
Observe each summand on the right in (\ref{hometa2}) has dimension 1.
These summands are mutually orthogonal so $m$ is equal to the 
dimension of
$E^*_2H_\eta $.
Now apparently $\mu_\eta $ is equal to the dimension of $E^*_2H_\eta 
$, as desired.
It follows the scalars (i), (ii) above are equal.
\hfill $\Box$\\
\bigskip


\section{Taut graphs and the local eigenvalues }  

\bigskip \noindent 
In this section we assume $\Gamma$ is bipartite with diameter $D 
\ge 4$,
valency $k \ge 3$, and eigenvalues $k=\theta_0 >\theta_1 > \cdots 
>\theta_D.$  Let $d = \lfloor D/2 \rfloor$.
We now recall the taut condition.
In \cite[Theorem 18]{curtin3} Curtin showed that $b_{2}(k-2) \ge (c_{2}-1)  \theta_{1}^{2}
$ with equality if and only if $\Gamma$ is 
2-homogeneous in the sense of Nomura 
\cite{Nom1}.  In \cite[Theorem 12]{curtin3} Curtin showed $\Delta \ge 0$, where 
$$
\Delta = (k-2)(c_{3}-1)-(c_{2}-1)p_{22}^{2}.
$$
In \cite[Lemma 3.8]{maclean1} MacLean proved
that 
\begin{equation}\label{bfb}
b_{3}\left(b_{2}(k-2)-(c_{2}-1)\theta_1^{2}\right)\left( 
b_{2} (k-2)-(c_{2}-1)\theta_d^{2}\right) \ge 
b_{1}\Delta(\theta_{1}^{2}-b_{2})(b_{2}-\theta_{d}^{2}).
\end{equation}
We remark that the inequality (\ref{bfb}) looks 
different from the inequality presented in \cite{maclean1}, but it is 
straightforward to show that these two inequalities are equivalent.
Recall from Lemma \ref{bipeigs} that $\theta_{1}^{2} > b_{2} > \theta_{d}^{2}$, 
so the last two factors on the right in (\ref{bfb}) are positive. 
Observe each factor in (\ref{bfb}) is nonnegative. 
From these comments we find that $\Gamma$ is
2-homogeneous if and only if $\Delta =0$ and equality holds in 
(\ref{bfb}).  
MacLean defined $\Gamma$ to be {\it taut} whenever $\Delta \not= 0$ 
and equality holds in (\ref{bfb}).

\bigskip
\noindent Fix $x \in X$. 
In this section we give a characterization of the taut 
condition in terms of the local
eigenvalues of $\Gamma$ with respect to $x$.
In order to motivate our results we sketch a proof of
(\ref{bfb}).

\bigskip
\noindent {\it Proof of} (\ref{bfb}). Fix $x \in X$ and let
$\eta_1, \eta_2,\ldots, \eta_{k_{2}}$
denote the corresponding local eigenvalues
of $\Gamma$.
Consider the sum
\begin{equation} \label{INTR0sum}
\sum_{i=1}^{k_{2}} (\eta_i-{\tilde \theta}_1)(\eta_i -{\tilde \theta}_d).
\end{equation}
We evaluate (\ref{INTR0sum}) in two ways.
First, by Theorem \ref{TYPEINEQ}
and since $\eta_1=p_{22}^{2}$, $\; \eta_{i}=b_{3}-1 \; (2 \le i \le k)$ we find
(\ref{INTR0sum}) is at most $(p_{22}^{2}-{\tilde \theta}_1)(p_{22}^{2} -{\tilde 
\theta}_d) + (k-1)(b_{3}-1-{\tilde \theta}_1)(b_{3}-1-{\tilde 
\theta}_d)$.
Second, we determine (\ref{INTR0sum}) by computing
$\sum_{i=1}^{k_{2}} \eta_i$ and $\sum_{i=1}^{k_{2}} \eta^2_i$. 
By Lemma \ref{multeqs}(ii),(iii) we find
$\sum_{i=1}^{k_{2}} \eta_i = 0 $ and
$\sum_{i=1}^{k_{2}} \eta^2_i = k_{2}p_{22}^{2}$.
From these comments the expression
(\ref{INTR0sum}) is equal to
$k_{2}(p_{22}^{2} +
{\tilde \theta}_1{\tilde \theta}_d) $.
We now have an inequality involving $k, k_{2}, p_{22}^{2}, b_{3},
{\tilde \theta}_1, {\tilde \theta}_d $; eliminating
${\tilde \theta}_1, {\tilde \theta}_d $ using (\ref{tildeth1d})
and simplifying using (\ref{KI}) and (\ref{p222}), we get
(\ref{bfb}).
\hfill $\Box$\\

\noindent In order to gain some insight into the case in which 
$\Gamma$ is taut, we examine the above proof. We begin with a 
definition.

\begin{definition}
\label{def:trelx}  \rm
With reference to
Definition \ref{A}, we say $\Gamma $ is
{\em spectrally taut with respect to
$x$} whenever $\eta_i$ is one of
${\tilde \theta}_1,
{\tilde \theta}_d $ for $k+1 \leq i \leq k_{2}$.
(The scalars $\eta_i$ and ${\tilde \theta}_i$ are from
Definition \ref{SUBGRAPH} and
Definition \ref{TYPE}, respectively.)
\end{definition}

 \noindent From the above proof of (\ref{bfb})
we routinely obtain the following.

\begin{theorem} 
\label{jk0char}
Let $\Gamma $ denote a bipartite distance-regular graph with
diameter $D\geq 4$ and valency $k \ge 3$. Then the following (i)--(iii) are equivalent.
\begin{description}
\item[(i)]
$\Gamma $ is taut.
\item[(ii)]
$\Delta \not= 0$ and $\Gamma $ is spectrally taut with
respect to each vertex.
\item[(iii)] $\Delta \not= 0$ and $\Gamma $ is spectrally taut with
respect to at least one vertex.
\end{description}
\end{theorem}

\begin{corollary}
With reference to Definition \ref{A}, assume that $\Gamma$ is taut.  Then
\begin{eqnarray}
\mbox{mult}_{{\tilde \theta}_{1}} &=& 
\frac{k(\theta_{1}^{2}-b_{2})(b_{2}(k-2)-(c_{2}-1)\theta_{d}^{2})}
{(\theta_{1}^{2}-\theta_{d}^{2})b_{2}c_{2}}, \label{mult1} \\
\mbox{mult}_{{\tilde \theta}_{d}} &=& 
\frac{k(\theta_{d}^{2}-b_{2})(b_{2}(k-2)-(c_{2}-1)\theta_{1}^{2})}
{(\theta_{d}^{2}-\theta_{1}^{2})b_{2}c_{2}}.  \label{multd}
\end{eqnarray}
Moreover each of $\mbox{mult}_{{\tilde \theta}_{1}}$, 
$\mbox{mult}_{{\tilde \theta}_{d}}$ is nonzero.  
\end{corollary}
{\it Proof.}  By Theorem \ref{jk0char}, $\Phi \subseteq \{{\tilde \theta}_1,
{\tilde \theta}_d \}$, where $\Phi$ is from Definition \ref{phidef}.
 Applying Lemma \ref{multeqs}(i),(ii), we find
$$ k + \mbox{mult}_{{\tilde \theta}_{1}} + \mbox{mult}_{{\tilde 
 \theta}_{d}} =k_{2}, \qquad p^{2}_{22} +(k-1)(b_{3}-1) + 
 {\tilde \theta}_{1}\mbox{mult}_{{\tilde \theta}_{1}} + 
 {\tilde \theta}_{d}\mbox{mult}_{{\tilde \theta}_{d}}=0.$$
Solving these two equations for $\mbox{mult}_{{\tilde \theta}_{1}}$, $\mbox{mult}_{{\tilde 
 \theta}_{d}}$ and simplifying using (\ref{KI}), (\ref{p222}) and 
(\ref{tildeth1d}), we obtain
(\ref{mult1}), (\ref{multd}).  Now we show $\mbox{mult}_{{\tilde \theta}_{1}}$, $\mbox{mult}_{{\tilde 
 \theta}_{d}}$ are nonzero.  Suppose $\mbox{mult}_{{\tilde 
 \theta}_{1}}=0.$
Since $\theta_{1}^{2}-b_{2} 
\not= 0$ by Lemma \ref{bipeigs}, by the form of 
(\ref{mult1}) we conclude that 
$b_{2}(k-2)-(c_{2}-1)\theta_{d}^{2}=0$.  Since equality holds in 
(\ref{bfb}), we have $\Delta =0$, contradicting Theorem \ref{jk0char}.
The proof that $\mbox{mult}_{{\tilde \theta}_{d}}$ is nonzero is 
similar.
\hfill $\Box$\\


\section{Taut graphs and the subconstituent algebra}

\noindent 
With reference to Definition \ref{A}, assume $D$ is odd. 
In this section we give two characterizations of the taut condition in 
terms of the subconstituent algebra $T$.  We begin with the following 
definition. 

\begin{definition}
\label{TAUTMOD}  \rm
With reference to Definition \ref{A}, assume $D$ is odd. 
We say $\Gamma $ is
{\it algebraically taut with respect to
$x$}
whenever every irreducible $T$-module with endpoint 2 is
thin with local eigenvalue
${\tilde \theta}_1$
or ${\tilde \theta}_d$.
 \end{definition}

\noindent The notions of spectrally taut and algebraically taut
are related as follows.

\begin{lemma}\label{UIUDTAUT}
With reference to Definition \ref{A}, assume $D$ is odd.  Then
the following (i), (ii) are equivalent.
\begin{description}
\item[(i)] $\Gamma $ is spectrally taut with respect to $x$.
\item[(ii)] $\Gamma $ is algebraically taut with respect to $x$.
\end{description}
\end{lemma}
{\it Proof.}
(i) $\Rightarrow $ (ii) Let $W$ denote an irreducible $T$-module
with endpoint 2. We show $W$ is thin with local eigenvalue
${\tilde \theta}_1$ or
${\tilde \theta}_d$.
Observe $E^*_2W$
is nonzero and invariant under $E^*_2A_{2}E^*_2$. Therefore
there exists a nonzero vector $v \in E^*_2W$ which is
an eigenvector for $E^*_2A_{2}E^*_2$. Let $\eta$ denote the
corresponding eigenvalue.  Let $Y$ denote the subspace of $V$ 
spanned by all the
irreducible $T$-modules with endpoint 1.  Observe $E^*_2W$ is orthogonal
to $V_{0} + Y$ so $v \in U$ by Definition \ref{U=V0+Y} . Now
$\eta$ is one of $\eta_{k+1}, \eta_{k+2}, \ldots, \eta_{k_{2}}$. By
this and Definition
\ref{def:trelx} we find $\eta$ is one of ${\tilde \theta}_1$,
${\tilde \theta}_d$.
By Theorem \ref{Gonew4}
we find $Mv$ is a thin irreducible $T$-module with endpoint
2 and local eigenvalue $\eta$. Observe $Mv \subseteq W$
so $Mv=W$ by the irreducibility of $W$.
Apparently
$W$ is thin with local eigenvalue $\eta$ and the result follows. \\
\noindent (ii) $\Rightarrow $ (i) Let $S$ denote the subspace of $V$ 
spanned by all the
irreducible $T$-modules with endpoint 2.  
Then \begin{equation}
S = H_{{\tilde \theta}_1}
+
H_{{\tilde \theta}_d} \qquad \qquad (\hbox{orthogonal direct sum}),
\label{SHH}
\end{equation}
where $H_{{\tilde \theta}_1}$ and
$H_{{\tilde \theta}_d}$
are from
Lemma \ref{T11}.  Let $Y$ denote the subspace of $V$ 
spanned by all the
irreducible $T$-modules with endpoint 1.
Recall $U$ denotes the orthogonal complement of $E^*_2V_{0} + E^*_2Y$ 
in 
$E^*_2V$.
Applying $E^*_2$ to each term in (\ref{SHH}), and evaluating the 
result using Lemma \ref{T11} and $E^*_2S=U$, we obtain 
\begin{equation}
U = U_{{\tilde \theta}_1} +
U_{{\tilde \theta}_d} \qquad \qquad
(\hbox{orthogonal direct sum}).
\label{UUU}
\end{equation}
Comparing (\ref{UUU}) and (\ref{Ubreakdown}) we find $\eta_i$ is one 
of ${\tilde \theta}_1$, ${\tilde \theta}_d$ for $k+1 \leq i \leq k_{2}$. 
Now $\Gamma$ is spectrally taut with respect to $x$ by Definition
\ref{def:trelx}.
\hfill $\Box$\\

\begin{definition}
\label{trx}  \rm
With reference to Definition \ref{A}, assume $D$ is odd.  We say $\Gamma$ is
{\it taut with respect to $x$} whenever the equivalent conditions 
(i), (ii) hold in Lemma
\ref{UIUDTAUT}.
\end{definition}

\noindent Combining
Theorem \ref{jk0char} and Definition \ref{trx}
we immediately obtain the following theorem.

\begin{theorem}\label{MODTIGHT}
Let $\Gamma$ denote a bipartite distance-regular graph with 
diameter $D \ge 4$ and valency $k \ge 3$.  Assume $D$ is odd.    
Then the following (i)--(iii) are equivalent.
\begin{description}
\item[(i)] $\Gamma$ is taut.
\item[(ii)] $\Delta \not= 0$ and $\Gamma$ is taut with respect to each 
vertex.
\item[(iii)] $\Delta \not= 0$ and $\Gamma$ is taut with respect to at 
least one vertex.
\end{description}
\end{theorem}

\noindent
For our final characterization of the taut condition, we will need 
the following two definitions.

\begin{definition}  \rm
Let $\Gamma$ denote a distance-regular graph with vertex set $X$ and 
diameter $D \ge 3$.  We say $\Gamma$ is an {\em antipodal 2-cover} whenever 
for all $x \in X$, there exists a unique vertex $y \in X$ such that 
$\partial(x,y)=D$.  In other words, $\Gamma$ is an antipodal 2-cover 
if and only if $k_{D}=1$.
\end{definition}

\begin{definition} \rm With reference to Definition \ref{A}, we say 
$\Gamma$ is {2-thin with respect to $x$} whenever every irreducible 
$T$-module with endpoint 2 is thin.
\end{definition}

\begin{theorem}
\label{thm:final}
Let $\Gamma$ denote a bipartite distance-regular graph with 
diameter $D \ge 4$ and valency $k \ge 3$.  Assume $D$ is odd.
Then the 
following (i)--(iii) are equivalent:
\begin{description}
\item[(i)] $\Gamma$ is taut or 2-homogeneous. 
\item[(ii)]  $\Gamma$ is an antipodal 2-cover and 2-thin with respect to each 
vertex.
\item[(iii)]  $\Gamma$ is an antipodal 2-cover and 2-thin with respect to at 
least one vertex.
\end{description}
\end{theorem}
{\it Proof.}  $(i)\Rightarrow(ii)$  Let $x$ denote a vertex of $\Gamma$.
First suppose $\Gamma$ is taut.  
Then $\Gamma$ is an antipodal 2-cover by \cite[Theorem 6.4]{maclean2}, and $\Gamma$ is 2-thin 
with respect to $x$ 
by Theorem \ref{MODTIGHT}.  Now suppose $\Gamma$ is 2-homogeneous.  
By \cite[Theorem 42]{curtin3}, $\Gamma$ is an antipodal 2-cover, and
$\Delta =0$ and equality holds in (\ref{bfb}).  Then $\Gamma$ is 
spectrally taut with respect to $x$ by the proof of (\ref{bfb}).  
Then $\Gamma$ is 2-thin with respect to $x$ by Definition 
\ref{TAUTMOD} and Lemma 
\ref{UIUDTAUT}.  \\
$(ii)\Rightarrow(iii)$  Clear.  \\
$(iii)\Rightarrow(i)$ By assumption there exists a vertex with 
respect to which $\Gamma$ is 2-thin.  Denote this vertex by $x$ and 
write $T=T(x)$.
Let $W$ denote an irreducible $T$-module of 
$\Gamma$ with endpoint 2.  Then the dimension of $W$ is $D-3$ by 
\cite[Lemma 14.1]{collins2}.  Now by Corollary \ref{dimj}, the local eigenvalue 
of $W$ is either ${\tilde \theta}_1$ or ${\tilde \theta}_d$.  Then 
$\Gamma$ is algebraically taut with respect to $x$ by Definition 
\ref{TAUTMOD}.  Thus $\Gamma$ is spectrally taut with respect to $x$ 
by Lemma \ref{UIUDTAUT}, so equality holds in (\ref{bfb}).  Thus 
$\Gamma$ is taut or 2-homogeneous.
\hfill $\Box$\\


\section{Directions for further research}

\noindent 
In this section we give some suggestions for further research.

\medskip
\noindent 
We start with a problem that we admit is quite general but
we believe is important.

\begin{problem}
\label{prob:gen}
\rm
With reference to Definition \ref{A}, assume that up to isomorphism 
there exist at most two irreducible $T$-modules with endpoint 2, and 
they are both thin. 
Investigate the combinatorial and algebraic 
implications of this assumption.
\end{problem}

\begin{remark}
\rm
With reference to Definition \ref{A},
assume $\Gamma$ is taut and $D$ is odd.
Then
$\Gamma$ satisfies the assumptions of Problem
\ref{prob:gen} by Theorem
\ref{MODTIGHT}.
\end{remark}

\begin{remark} 
\rm
With reference to Definition \ref{A},
assume $\Gamma$ is $Q$-polynomial. 
Then
$\Gamma$ satisfies the assumptions of Problem
\ref{prob:gen}
\cite[Section 14]{caugh2}.
\end{remark}

\begin{remark}
\label{rem:rec}
\rm
Assume $\Gamma$ satisfies the assumptions of
Problem 
\ref{prob:gen}.
Then using
\cite[Theorem 13.1]{curtin2} we can recursively obtain 
the intersection
numbers of $\Gamma$ in terms of the diameter $D$,
the local
eigenvalues of the 
$T$-modules mentioned in Problem
\ref{prob:gen}, and the multiplicities with which
these modules appear in $V$. The resulting formulae
are not attractive however.
\end{remark}

\begin{problem}
\rm
Assume $\Gamma$ satisfies the assumptions of
Problem 
\ref{prob:gen}.
Obtain the intersection numbers of
$\Gamma$ in closed form 
as attractive rational expressions involving
$D$ and at most four complex parameters. 
\end{problem}


\begin{problem}
\rm
With reference to Definition \ref{A}, let $W$ denote a thin 
irreducible $T$-module with endpoint 2. Let 
$\eta$  denote the local eigenvalue of $W$ and assume
${\tilde \theta}_1 < \eta < {\tilde \theta}_d$.
Describe the 
structure of $W$ along the lines of Section 13 in the
present paper.
\end{problem}

\noindent The following problem is of interest in view of
Theorem
\ref{jk0char}.

\begin{problem}
\rm
With reference to 
Definition \ref{A} and Definition
\ref{U=V0+Y},
assume $D$ is 
even and there exists
a nonzero vector 
$v \in U$ that is an eigenvector for 
$E_{2}^{*}A_{2}E_{2}^{*}$ with  
eigenvalue $b_3-1$.
Show the $T$-module $Tv$ is irreducible with basis
\begin{eqnarray*}
\lbrace E^*_{i+2}A_iv \,|\,0 \leq i \leq D-3\rbrace
\quad \cup \quad 
\lbrace E^*_iA_{i+2}v \,|\,3 \leq i \leq D-3,\;\; i\; \mbox{odd}\rbrace.
\end{eqnarray*}
Compute the matrix representing $A$ with respect to this basis.
\end{problem} 

\begin{problem}
\label{def:pstaut}
\rm
Let $\Gamma$ denote a bipartite distance-regular
graph with diameter $D\geq 4$.
Let $E,F$ denote  pseudo primitive idempotents
of $\Gamma$ \cite{tw}. We say the pair
$E,F$ is 
{\it taut}
whenever $E\circ F$ is a linear combination
of at most two pseudo primitive idempotents of $\Gamma$.
Find all the taut pairs of pseudo primitive idempotents
of $\Gamma$. See 
\cite{pascter} for related results concerning tight
pairs of pseudo primitive idempotents.
\end{problem}

\begin{conjecture}
\label{con:pseudo}
\rm
With reference to Definition \ref{A}, assume that up to isomorphism 
there exist exactly two irreducible $T$-modules with endpoint 2, and 
they are both thin. 
Pick $\xi,\chi \in \C \cup \infty$ such
that $\tilde \xi$ and
$\tilde \chi$ are the local eigenvalues for these $T$-modules.
Let $E$ (resp. $F$) denote 
a pseudo primitive idempotent \cite{tw}
of $\Gamma$ for
$\xi$ (resp. $\chi$). 
Then the pair $E,F$ is taut in
the sense of Problem 
\ref{def:pstaut}.
\end{conjecture}

\begin{conjecture}
\rm
Let $\Gamma$ denote a taut distance-regular graph
with odd diameter $D\geq 5$.
Recall $\Gamma$ is an antipodal 2-cover by
Theorem 
\ref{thm:final}; we conjecture that the antipodal
quotient of $\Gamma$ is $Q$-polynomial. 
\end{conjecture}

\begin{problem}
\label{prob:dep}
\rm
Let $\Gamma=(X,R)$ denote a bipartite
distance-regular graph with diameter $D\geq 4$
and eigenvalues $\theta_0 > \theta_1 > \cdots > \theta_D$.
Fix
$x,y \in X$ 
at distance $\partial (x,y)=2$.
For $0 \leq i ,j\leq D$ define
$w_{ij}=\sum {\hat z}$,
where the sum is over all  $z \in X$
such that $\partial(x,z)=i$ and $\partial(y,z)=j$.
Show that the following are equivalent:
(i) $\Gamma$ is taut;
(ii) The vectors $E_1{\hat x}$, 
$E_1{\hat y}$, 
$E_1w_{11}$, 
$E_1w_{22}$ are linearly dependent and
the vectors $E_d{\hat x}$, 
$E_d{\hat y}$, 
$E_dw_{11}$, 
$E_dw_{22}$ are linearly dependent,
where 
$d = \lfloor D/2 \rfloor$.
\end{problem}

\begin{problem}
\rm
Let $\Gamma=(X,R)$ denote a taut distance-regular
graph with odd diameter $D\geq 5$.
Fix 
$x,y \in X$ 
at distance $\partial (x,y)=2$,
and let the vectors 
$w_{ij}$ be as in 
Problem
\ref{prob:dep}.
For $2 \leq i \leq D-2$ define
$w^+_{ii}=\sum \vert \lbrace r \in X\;|\; \partial (r,x)=1,
\partial(r,y)=1,\partial(r,z)=i-1\rbrace \vert
\,\hat{z}$,
where the sum is over all $z \in X$ such that
$\partial (x,z)=i$ and $\partial (y,z)=i$.
Show that the vectors
\begin{eqnarray*}
\lbrace w_{ij} \;|\;0 \leq i ,j\leq D,\; |i-j|=2\rbrace
\quad \cup \quad 
\lbrace w_{ii} \;|\;1\leq i \leq D-1\rbrace
\quad \cup \quad 
\lbrace w^+_{ii} \;|\;2\leq i \leq D-2\rbrace
\end{eqnarray*}
form a basis for an $A$-invariant subspace of
the standard module. Find the matrix representing
$A$ with respect to this basis.
\end{problem}



\end{document}